\documentclass{article}

\usepackage{amsfonts}

\usepackage{amssymb}

\newcommand{\marcio}[1]{\stackrel{#1}{_{\longmapsto}}}

\title{Flow: the Axiom of Choice is independent from the Partition Principle}

\author{Adonai S. Sant'Anna \and Ot\'avio Bueno \and Marcio P. P. de Fran\c ca \and Renato Brodzinski}

\date{\small For correspondence: adonaisantanna@gmail.com}

\begin{document}

\newtheorem{definicao}{Definition}
\newtheorem{teorema}{Theorem}
\newtheorem{lema}{Lemma}
\newtheorem{corolario}{Corollary}
\newtheorem{proposicao}{Proposition}
\newtheorem{axioma}{Axiom}
\newtheorem{observacao}{Observation}




\maketitle

\begin{abstract}

We introduce a general theory of functions called Flow. We prove ZF, non-well founded ZF and ZFC can be immersed within Flow as a natural consequence from our framework.  The existence of strongly inaccessible cardinals is entailed from our axioms. And our first important application is the introduction of a model of Zermelo-Fraenkel set theory where the Partition Principle (PP) holds but not the Axiom of Choice (AC). So, Flow allows us to answer to the oldest open problem in set theory: if PP entails AC.\\

\begin{center}
{\sc This is a full preprint to stimulate criticisms before we submit a final version for publication in a peer-reviewed journal.}
\end{center}
\end{abstract}


\section{Introduction}

It is rather difficult to  determine when functions were born, since mathematics itself changes along history. Specially nowadays we find different formal concepts associated to the label {\em function\/}. But an educated guess could point towards Sharaf al-D\={\i}n al-T\={u}s\={\i} who, in the 12th century, not only introduced a `dynamical' concept which could be interpreted as some notion of function, but also studied how to determine the maxima of such functions \cite{Hogenduk-89}.

All usual mathematical approaches for physical theories are based on either differential equations or systems of differential equations whose solutions (when they exist) are either functions or classes of functions (see, e.g., \cite{Sobolev-11}). In pure mathematics the situation is no different. Continuous functions, linear transformations, homomorphisms, and homeomorphisms, for example, play a fundamental role in topology, linear algebra, group theory, and differential geometry, respectively. Category theory \cite{MacLane-94} emphasizes such a role in a very clear, elegant, and comprehensive way. As remarked by Marquis \cite{Marquis-18}, category theory allows us to distinguish between canonical and noncanonical maps, in a way which is not usually achieved within purely extensional set theories, like ZFC. And canonical maps ``constitute the highway system of mathematical concepts''. Concepts of symmetry are essential in pure and applied mathematics, and they are stated by means of group transformations \cite{Sternberg-95}.

Some authors suggest that functions are supposed to play a strategic role into the foundations of mathematics \cite{vonNeumann-25} and even mathematics teaching \cite{Klein-16}, rather than sets. The irony of such perspective lies in the historical roots of set theories. Georg Cantor's seminal works on sets were strongly motivated by Bernard Bolzano's manuscripts on infinite multitudes called {\em Menge\/} \cite{Trlifajova-18}. Those collections were supposed to be conceived in a way such that the arrangement of their components is unimportant. However, Bolzano insisted on an Euclidian view that the whole should be greater than a part, while Cantor proposed a quite different approach: to compare infinite quantities we should consider a {\em one-to-one correspondence\/} between collections. Cantor's concept of collection (in his famous {\em Mengenlehre\/}) was strongly committed to the idea of function. Subsequent formalizations of Cantor's ``theory'' were developed in a way such that all strategic terms were associated to an intended interpretation of collection. The result of that effort is a strange phenomenon from the point of view of theories of definition: Padoa's principle allows us to show domains of functions are definable from functions themselves \cite{daCosta-01} \cite{daCosta-02} \cite{Sant'Anna-14}. Thus, even in extensional set theories functions seem to play a more fundamental role than sets.

Sets can be intuitively viewed as the result of a process of collecting objects. An object is collected if it is assigned to a given set. But the fundamental mechanism here is {\em to attribute\/} something to a certain collection, a process which resembles the role a function is supposed to perform. From another point of view, we should recall that sets and functions are meant to correspond to an intuitive notion of properties. Usually properties allow to define either classes or sets (like the Separation Schema in ZFC). But another possibility is that properties correspond to functions. Talking about objects that have a given property $P$ corresponds to associate certain objects to a label which represents $P$; and any other remaining objects are supposed to be associated to a different label. The correspondence itself between $P$ and a given label does have a functional, rather than a set-theoretical, appeal. And usually, those labels are called sets. So, why do we need sets? Why can't we deal primarily with functions? In other words, why can't we label those intended properties with functions instead of sets?

Another issue with extensional set theories like ZFC and NBG is the fact that usually sets are solely `generated' from a hierarchical procedure characterized by postulates like `pair', `power', `replacement', and `union'. Few exceptions are axioms like `choice' and `separation'. But even those last cases depend on the former ones to work. That character of standard set theories generates a lot of discussions regarding independent formulas like those which guarantee the existence of (even weakly) inaccessible cardinals. Thus, if someone intends to make category theory work within a ZF-style set theory, it is necessary to postulate an independent formula which guarantees the existence of some sort of universe, like Grothendieck's \cite{MacLane-69} or hierarchies of constructible sets \cite{Jech-03}.

It could be thought that category theory itself provides a framework to develop this sort of functional approach we intend to investigate. Category theory deals primarily with morphisms \cite{MacLane-94} which are commonly associated to some notion of function (at least in a vast amount of applications). For example, the category of topological spaces, whose objects are topological spaces, may have their morphisms described as homotopy classes of continuous maps, and not functions. But even in that case the concept of function (map) is somehow present.

However, even morphisms, functors, and natural transformations have domains and codomains, and so we still wouldn't have the appropriate framework to develop what we have in mind. That happens because Category Theory seems to be somehow historically committed to some sort of set-theoretic intuition about what a function is supposed to be: an `agent' which is supposed to `act' on objects which inhabit some {\em stage\/} commonly understood as a collection. Such a stage is usually referred to as a domain. And we want to investigate the possibility of getting rid of domains. In other words, we are pursuing a purely functional formal framework where collections truly play a secondary role.

String diagrams \cite{Dixon-13} (whose edges and vertices can eventually be interpreted as morphisms in a monoidal category) introduce some ideas which seem to intersect with our own. A remarkable feature of string diagrams is that edges need not be connected to vertices at both ends. More than that, unconnected ends can be interpreted as inputs and outputs of a string diagram (with important applications in computer science). But within our framework {\em no\/} function has any domain whatsoever. Besides, the usual way to cope with string diagrams is by means of a discrete and finitary framework, while in Flow (our framework) we have no need for such restrictions.

A more radical proposal for mathematical foundations where functions play a fundamental role is the Theory of Autocategories \cite{Guitart-14}. Autocategories are developed with the aid of autographs, where arrows (which work as morphisms) are drawn between arrows with no need for objects. And once more we see a powerful idea concerning functions is naturally emerging in different places nowadays.

This paper is motivated by von Neumann's original ideas \cite{vonNeumann-25} and a variation of them \cite{Sant'Anna-14}, and related papers as well (\cite{daCosta-01} \cite{daCosta-02}). In \cite{Sant'Anna-14} it was provided a reformulation of von Neumann's `functions' theory (termed $\mathcal N$ theory). Nevertheless, we develop here a whole new approach. We close this Introduction with our main ideas from an intuitive point of view, including our main contribution. Locutions like ``postulate'' and ``axiom'' are used interchangeably, as synonyms.

Flow is a first order theory with equality, where the intended interpretation of its terms is that of a monadic function. All terms are called {\em functions\/}. Our framework has one primitive concept (besides equality): a functional letter $f_1^2$ where $f_1^2(f,x) = y$ is abbreviated as $y = f(x)$. We call $y$ the {\em image\/} of $x$ by $f$.

Every $f$ has a unique image $f(t)$ for any $t$. Thus, our functions have no domain whatsoever, although we are able to define concepts like domain and range. Besides, there are two special functions: $\underline{\mathfrak{0}}$ and $\underline{\mathfrak{1}}$. For any $t$ we have $\underline{\mathfrak{0}}(t) = \underline{\mathfrak{0}}$ and $\underline{\mathfrak{1}}(t) = t$. So, for all that matters, $\underline{\mathfrak{1}}$ behaves like an identity function and $\underline{\mathfrak{0}}$ works as a null function.

A {\em Self-Reference Axiom\/} says $f(f) = f$, for any $f$. Any functions $f$ and $g$ which share the same images $f(t)$ and $g(t)$, for any $t$, are the same. That is provable as a consequence from a {\em Weak Extensionality Axiom\/} and our self-reference postulate.

If $t\neq f$ and $f(t)\neq\underline{\mathfrak{0}}$, we say $f$ {\em acts on\/} $t$ and denote this by $f[t]$. No $f$ acts on itself, according to the self-reference axiom. We use this concept for defining a membership relation: $t\in f$ iff $f[t]$ and some conditions are imposed for both $f$ and $t$. In this sense, our self-reference axiom works as very weak form of regularity, since it entails there is no $f$ such that $f\in f$.

For {\em any\/} $f$ and $g$ there is a binary {\em $\mathfrak{F}$-composition\/} $f\circ g$ (which is associative only on a portion of the universe of discourse of Flow), such that $(f\circ g)(t) = f(g(t))$, except when $t = f \circ g$, $t = f$ and $t = g$. If $f\circ g\neq f$, then $(f\circ g)(f) = \underline{\mathfrak{0}}$; and the same happens to $f\circ g \neq g$. Thus, $\mathfrak{F}$-composition does not rely on any consideration regarding domains or codomains of the functions to be $\mathfrak{F}$-composed. Due to the self-reference postulate, $(f\circ g)(f\circ g) = f\circ g$.

The {\em $\mathfrak{F}$-successor\/} $\sigma(f)$ (abbreviated as $\sigma_f$) of an $f$ is supposed to be a function $g$ which shares the same images of $f$ for any $t\neq g$, in a way such that $f\neq g$. If there is no $g$ like that, then $\sigma_f$ is simply $\underline{\mathfrak{0}}$, except for the cases when $f$ is either $\sigma$ itself (due to self-reference) or $\underline{\mathfrak{0}}$ (since there is a theorem which states $f(\underline{\mathfrak{0}}) = \underline{\mathfrak{0}}$ for any $f$). In those cases, $\sigma_\sigma = \sigma$ and $\sigma_{\underline{\mathfrak{0}}} = \underline{\mathfrak{0}}$, where $\sigma$ itself is a function.

There is a function $\varphi_0$ which is identifiable to the empty set of extensional set theories in a very precise manner. Its existence is granted by $\mathfrak{F}$-composition. $\varphi_0$ and $\underline{\mathfrak{0}}$ are the only functions who do not act on any term, although $\varphi_0\neq\underline{\mathfrak{0}}$. Besides, $\sigma_{\varphi_0} = \varphi_1$, where $\varphi_1$ is identifiable to the unitary set whose only member is the empty set; $\sigma_{\varphi_1} = \varphi_2$; and so on. Those $\varphi_n$ are our first ordinals.

A {\em restriction\/} $g$ of a given $f$ is such that for any $t$, if $g[t]$, then $f[t]$ and $g(t) = f(t)$. Restrictions may be definable (but not always) from an axiom which resembles the Separation Scheme of ZFC. Many restrictions $g$ of a given function $f$ are supposed to be determined under the influence of a formula $F$: the values $g(t)$ depend on the values $f(t)$ and whether or not we have $F(t)$. The existence of such restrictions depends on considerations regarding two different kinds of functions: those who act on `many' terms and those who do not. Functions who act on `many' terms are called {\em comprehensive\/} (some of them correspond to the usual notion of a {\em proper class\/}). A function $f$ is comprehensive iff there is a restriction $g$ of $f$ such that $\sigma_g = \underline{\mathfrak{0}}$ and $g\neq\underline{\mathfrak{0}}$. The remaining functions are called {\em uncomprehensive\/}. We denote $g$ as $f\big|_F$. When $f$ is uncomprehensive, there are no constraints whatsoever for determining $g$. But when $f$ is comprehensive, we should be very careful. In that case, we demand formula $F$ be such that, for any $t$, $F(t)$ entails $t$ is uncomprehensive. But even with such a constraint, we are still able to define uncomprehensive functions which act on comprehensive terms. That allows us to define a Grothendieck Universe.

A term $g$ {\em lurks\/} $f$ iff: $\forall x(g[x]\Rightarrow f[x])$, or $x=f(t)$ for some $t$ {\em and\/} either $f[g(x)]$ or $g(x)=f(t)$ for some $t$. That concept generalizes the usual notion of subset.

The {\em full power\/} of any uncomprehensive function $f$ is the uncomprehensive function $\mathfrak{p}(f)$ which acts on every possible $g$ that lurks $f$. For example, if $f$ acts on just two distinct terms $x$ and $y$, such that $f(x) = x$ and $f(y) = y$, then $\mathfrak{p}(f)$ acts on the next functions: $f_1$, $f_2$, $f_3$, $f_4$, $f_5$, $f_6$, $f_7$, $f_8$, and $f_9$, where: $f_1(x) = x$, $f_1(y) = y$; $f_2(x) = x$, $f_2(y) = x$; $f_3(x) = y$, $f_3(y) = y$; $f_4(x) = y$, $f_4(y) = x$; $f_5(x) = y$, $f_5(y) = \underline{\mathfrak{0}}$; $f_6(y) = x$, $f_6(x) = \underline{\mathfrak{0}}$; $f_7(x) = x$, $f_7(y) = \underline{\mathfrak{0}}$; $f_8(y) = y$, $f_8(x) = \underline{\mathfrak{0}}$; $f_9 = \varphi_0$. Thus, $\mathfrak{p}(f)$ unavoidably acts on $f$ (here represented by $f_1$) and $\varphi_0$. We show this is a significant generalization of the concept of power set. Both full power and the usual notion of power set (called here {\em restricted power\/}) are definable within our framework. Actually, that is one of the methodological and epistemological advantages of Flow for solving the problem of Partition Principle (PP) {\em versus\/} Axiom of Choice (AC). More than that, Flow allows us to introduce models of ZF with AC or its negation, with the Axiom of Regularity or its negation, among other possibilities. On the other hand, our framework is not supposed to be simply one more technique for creating models of ZF besides those already stated in the literature \cite{Abian-78} \cite{Devlin-93} \cite{Jech-03} \cite{Kanamori-03} \cite{Steel-07}. Our main purpose is to introduce and develop this general theory of functions with merits by its own.

{\em Creation Axiom\/} helps to build functions which are no restrictions of $\underline{\mathfrak{1}}$. Although we do not need the concepts of domain and range of a function, those are easily definable. Thus we are able to translate the concepts of injective and surjective functions within our framework.

Finally, a $\mathfrak{F}$-Choice Postulate grants the existence of a quite useful injection $g$ that lurks any surjection $f$. That postulate entails PP but not AC. Thus, Flow is used to build a model where both PP and the negation of AC hold.

When Ernst Zermelo introduced AC, his motivation was the Partition Principle \cite{Banaschewski-90}. Indeed, AC entails PP. But since 1904 (when Zermelo introduced AC) it is unknown whether PP implies AC. Our answer to that question is negative, thanks to Theorems \ref{teoremadeparticao} and \ref{onzedeagostode2020}. Our point is that ZF, ZFC, all their variants, and almost all their respective models (with a few exceptions like those in \cite{Abian-78} and \cite{Devlin-93}) are somehow committed to a methodological and epistemological character which forces us to see sets as collections of some sort. Within Flow that does not happen, since our framework drives us to see sets as special cases of functions. All sets in Flow are restrictions of $\underline{\mathfrak{1}}$. Nevertheless, functions who are no restrictions of $\underline{\mathfrak{1}}$ play an important role with consequences over sets. An analogous situation takes place with the well known Reflexion Principle in model theory \cite{Steel-07}: some features of the von Neumann Universe motivate mathematicians to look for new axioms which grant the existence of strongly inaccessible cardinals. Flow, however, provides a new way for coping with the metamathematics of ZF.

\section{Flow theory}

Flow is a first-order theory $\mbox{\boldmath{$\mathfrak{F}$}}$ with identity \cite{Mendelson-97}, where $x=y$ reads ``$x$ is equal to $y$'', and $\neg(x = y)$ is abbreviated as $x\neq y$. Flow has one functional letter $f_1^2$ (termed {\em evaluation\/}), where $f_1^2(f,x)$ is a term, if $f$ and $x$ are terms. If $y = f_1^2(f,x)$, we abbreviate this by $f(x) = y$, and read ``$y$ is the {\em image\/} of $x$ by $f$''. All terms are called {\em functions\/}. Such a terminology seems to be adequate under the light of our intended interpretation: functions are supposed to be terms which `transform' terms into other terms. Since we are assuming identity, our functions cannot play the role of non-trivial relations. Suppose, for example, $f(x) = y$ and $f(x) = y'$, which are abbreviations for $f_1^2(f,x) = y$ and $f_1^2(f,x) = y'$, respectively. From transitivity of identity, we have $y = y'$. So, for any $f$ and any $x$, there is one single $y$ such that $f(x) = y$.

Lowercase Latin and Greek letters denote functions, with the sole exception of two specific functions to appear in the next pages, namely, $\underline{\mathfrak{0}}$ and $\underline{\mathfrak{1}}$. Uppercase Latin letters are used to denote formulas or predicates (which are eventually defined). Any explicit definition in Flow is an abbreviative one, in the sense that for a given formula $F$, the {\em definiendum\/} is a metalinguistic abbreviation for the {\em definiens\/} given by $F$. Eventually we use bounded quantifiers. If $P$ is a predicate defined by a formula $F$, we abbreviate $\forall x(P(x)\Rightarrow G(x))$ and $\exists x(P(x)\wedge G(x))$ as $\forall_P x(G(x))$ and $\exists_P x(G(x))$, respectively; where $G$ is a formula. Finally, $\exists!x(G(x))$ is a metalinguistic abbreviation for $\exists x\forall y(G(y)\Leftrightarrow y = x)$.

Postulates {\sc F1}$\sim${\sc F11} of $\mbox{\boldmath{$\mathfrak{F}$}}$ are as follows.

\begin{description}

\item[\sc F1 - Weak Extensionality] $\forall f \forall g (((f(g) = f \wedge g(f) = g) \vee (f(g) = g \wedge g(f) = f)) \Rightarrow f = g))$.

\end{description}

If $f(g) = f$ we say $f$ is {\em rigid with\/} $g$. If $f(g) = g$ we say $f$ is {\em flexible with\/} $g$. If both $f$ and $g$ are rigid (flexible) with each other, then $f = g$.

\begin{description}

\item[\sc F2 - Self-Reference] $\forall f (f(f) = f)$.

\end{description}

Any $f$ is flexible {\em and\/} rigid with itself. That may sound a strong limitation. But that feature is quite useful for our purposes.

\begin{observacao}\label{Russell}

Let $y$ be a function such that $\forall x(y(x) = r \Leftrightarrow x(x)\neq r)$. What about $y(y)$? If $x = y$, then $y(y) = r \Leftrightarrow y(y)\neq r$ (Russell's paradox). But that entails $y(y)\neq y(y)$, another contradiction. However, our Self-Reference postulate does not allow us to define such an $y$, since $x(x) = x$ for any $x$. That is not the only way to avoid such an inconsistency. Our Restriction axiom (some pages below) does the same work. That opens the possibility of variations of Flow where Self-Reference does not hold.
\end{observacao}

Our first theorem states any $f$ can be identified by its images $f(x)$.

\begin{teorema}\label{igualdadefuncoes}
$\forall f\forall g(f = g\Leftrightarrow\forall x (f(x) = g(x)))$.
\end{teorema}

\begin{description}
\item[Proof] From substitutivity of identity in $f(x) = f(x)$, proof of the $\Rightarrow$ part is straightforward: if $f = g$, then $f(x) = g(x)$, for any $x$. For the $\Leftarrow$ part, suppose, for any $x$, $f(x) = g(x)$. If $x = f$, $f(f) = g(f)$; if $x = g$, $f(g) = g(g)$. From {\sc F2}, $f(f) = f$ and $g(g) = g$. So, $g(f) = f$ and $f(g) = g$. {\sc F1} entails $f = g$.
\end{description}

Axioms {\sc F1} and {\sc F2} could be rewritten as one single formula:

\begin{description}

\item[\sc F1' - Alternative Weak Extensionality] $\forall f \forall g (((f(g) = f \wedge g(f) = g) \vee (f(g) = g \wedge g(f) = f)) \Leftrightarrow f = g))$.

\end{description}

{\sc F2} is a consequence from {\sc F1'}. Ultimately, $f=g$ entails $f(g) = f$ (from {\sc F1'}). And substitutivity of identity entails $f(f) = f$. But we prefer to keep axioms {\sc F1} and {\sc F2} (instead of {\sc F1'}) to smooth away some discussions. From {\sc F1} and {\sc F2}, we can analogously see that {\sc F1'} is a nontrivial theorem.

\begin{description}

\item[\sc F3 - Identity] $\exists f \forall x (f(x) = x)$.

\end{description}

There is at least one $f$ such that, for any $x$, we have $f(x) = x$. Any function $f$ which satisfies {\sc F3} is said to be an {\em identity function\/}.

\begin{teorema}\label{umbarraunico}
Identity function is unique.
\end{teorema}

\begin{description}
\item[Proof] Suppose both $f$ and $g$ satisfy {\sc F3}. Then, for any $x$, $f(x) = x$ and $g(x) = x$. Thus, $f(g) = g$ and $g(f) = f$. Hence, according to {\sc F1}, $f = g$.
\end{description}

In other words, there is one single $f$ which is flexible to every term. In that case we simply say $f$ is {\em flexible\/}. That means ``flexible'' and ``identity'' are synonyms.

\begin{description}

\item[\sc F4 - Rigidness] $\exists f \forall x (f(x) = f)$.

\end{description}

There is at least one $f$ which is rigid with any function. Observe the symmetry between {\sc F3} and {\sc F4}! Any $f$ which satisfies last postulate is said to be {\em rigid\/}.

\begin{teorema}\label{unicidadezerobarra}
The rigid function is unique.
\end{teorema}

\begin{description}
\item[Proof] Let $f$ and $g$ satisfy axiom {\sc F4}. Then, for any $x$, $f(x) = f$ and $g(x) = g$. Thus, $f(g(x)) = f(g) = f$ and $g(f(x)) = g(f) = g$. From {\sc F1}, $f = g$.
\end{description}

Now we are able to introduce new terminology. Term $\underline{\mathfrak{1}}$ is the identity (flexible) function, while $\underline{\mathfrak{0}}$ is the rigid function, since we proved they are both unique. So,
$$\forall x (\underline{\mathfrak{1}}(x) = x\;\;\wedge\;\;\underline{\mathfrak{0}}(x) = \underline{\mathfrak{0}}).$$
Since $f(x) = y$ says $f_1^2(f,x) = y$, {\sc F3} states there is an $f$ such that, for any $x$, $f_1^2(f,x) = x$, while {\sc F4} says there is $f$ where $f_1^2(f,x) = f$. If we do not grant the existence of other functions, it seems rather difficult to prove $\underline{\mathfrak{0}}\neq\underline{\mathfrak{1}}$.

\begin{teorema}\label{teoremamarcio1}
$\underline{\mathfrak{0}}$ is the only function which is rigid with $\underline{\mathfrak{0}}$.
\end{teorema}

\begin{description}
\item[Proof] This theorem says $\forall x (x\neq \underline{\mathfrak{0}} \Rightarrow x(\underline{\mathfrak{0}}) \neq x)$, i.e., $\forall x (x(\underline{\mathfrak{0}}) = x \Rightarrow x = \underline{\mathfrak{0}})$. But $\underline{\mathfrak{0}}(x) = \underline{\mathfrak{0}}$. So, if  $x(\underline{\mathfrak{0}}) = x \wedge \underline{\mathfrak{0}}(x) = \underline{\mathfrak{0}}$, then $x = \underline{\mathfrak{0}}$ ({\sc F1}).
\end{description}

\begin{teorema}\label{teoremamarcio2}
$\underline{\mathfrak{1}}$ is the only function which is flexible with $\underline{\mathfrak{1}}$.
\end{teorema}

\begin{description}
\item[Proof] This theorem says $\forall x (x\neq \underline{\mathfrak{1}} \Rightarrow x(\underline{\mathfrak{1}}) \neq \underline{\mathfrak{1}})$, i.e., $\forall x (x(\underline{\mathfrak{1}}) = \underline{\mathfrak{1}} \Rightarrow x = \underline{\mathfrak{1}})$. But $\underline{\mathfrak{1}}(x) = x$. So, if  $x(\underline{\mathfrak{1}}) = \underline{\mathfrak{1}} \wedge \underline{\mathfrak{1}}(x) = x$, then $x = \underline{\mathfrak{1}}$ ({\sc F1}).
\end{description}

For more details about $x(\underline{\mathfrak{0}})$ and $x(\underline{\mathfrak{1}})$, for any $x$, see Theorems \ref{fechamentorenato} and \ref{composicaodexcomumbarra}.

\begin{definicao}\label{composicao}
Given $f$ and $g$, the {\em $\mathfrak{F}$-composition\/} $h = f\circ g$, if it exists, must satisfy the next conditions:

{\sc (i)} $h\neq\underline{\mathfrak{0}}$; {\sc (ii)} $\forall x((x\neq f \wedge x\neq g \wedge x\neq h)\Rightarrow h(x) = f(g(x)))$; \\{\sc (iii)} $h\neq f \Rightarrow h(f) = \underline{\mathfrak{0}}$; {\sc (iv)}  $h\neq g \Rightarrow h(g) = \underline{\mathfrak{0}}$; \\{\sc (v)} $(g\neq h\wedge f\neq h \wedge g\neq\underline{\mathfrak{1}}\wedge f\neq\underline{\mathfrak{1}})\Rightarrow (f(g(h)) = \underline{\mathfrak{0}}\vee g(h) = \underline{\mathfrak{0}})$

\end{definicao}

So, {\sc (i)} No $\mathfrak{F}$-composition is $\underline{\mathfrak{0}}$. {\sc (ii)} $\mathfrak{F}$-composition behaves like standard notions of composition between functions up to self-reference. Nevertheless, contrary to usual practice in standard set theories, we can define $\mathfrak{F}$-composition between any two functions. {\sc (iii and iv)} No $\mathfrak{F}$-composition acts on any of its factors. Finally, item {\sc (v)} is a technical constraint to avoid ambiguities in the calculation of $f\circ g$. Either $g$ does not act on the $\mathfrak{F}$-composition $h$ or, if it acts, then $f$ does not act on $g(h)$.

\begin{description}

\item[\sc F5 - $\mathfrak{F}$-Composition] $\forall f \forall g\exists !h (h = f\circ g)$.

\end{description}

We never calculate $(f\circ g)(f\circ g)$ as $f(g(f\circ g))$, since $(f\circ g)(f\circ g)$ is $f\circ g$, according to {\sc F2}. Given $f$ and $g$, we can `build' $f\circ g$ through a three-step process: {\sc (i)} First we establish a label $h$ for $f\circ g$; {\sc (ii)} Next we evaluate $f(g(x))$ for any $x$ which is different of $f$ and $g$. By doing that we are assuming those $x$ are different of $h$. Thus, if such a choice of $x$ entails $f(g(x)) = y$ for a given $y$, then $h(x) = y$; {\sc (iii)} Next we evaluate the following possibilities: is $h$ equal to either $f$, $g$ or something else? If $h\neq g$, then $h(g)$ is supposed to be $\underline{\mathfrak{0}}$. If $h(g) = \underline{\mathfrak{0}}$ entails a contradiction, then $h$ is simply $g$. An analogous method is used for assessing if $h$ is $f$. Eventually, $h$ is neither $f$ nor $g$, as we can see in the next theorems.

\begin{teorema}\label{primeirophizero}

There is a unique $h$ such that $h\neq\underline{\mathfrak{0}}$ but $h(x) = \underline{\mathfrak{0}}$ for any $x\neq h$.

\end{teorema}

\begin{description}
\item[Proof] From Definition \ref{composicao} and {\sc F5}, $\underline{\mathfrak{0}}\circ\underline{\mathfrak{0}}$ is a unique $h\neq\underline{\mathfrak{0}}$ such that, for any $x$ where $x\neq \underline{\mathfrak{0}}$ and $x\neq h$, $h(x) = \underline{\mathfrak{0}}(\underline{\mathfrak{0}}(x))=  \underline{\mathfrak{0}}(\underline{\mathfrak{0}}) = \underline{\mathfrak{0}}$. Since $h\neq\underline{\mathfrak{0}}$, then {\sc F5} entails $h(\underline{\mathfrak{0}}) = \underline{\mathfrak{0}}$. Thus, $h(x)$ is $\underline{\mathfrak{0}}$ for any $x\neq h$, while $h$ itself is different of $\underline{\mathfrak{0}}$.
\end{description}

Such $h$ of last theorem (which does not conflict neither with Theorem \ref{unicidadezerobarra} nor with Theorem \ref{teoremamarcio1}) is labeled with a special symbol, namely, $\varphi_0$. So, $\underline{\mathfrak{0}}\circ\underline{\mathfrak{0}} = \varphi_0$, where $\varphi_0\neq\underline{\mathfrak{0}}$. If the reader is intrigued by the subscript $_0$ in $\varphi_0$, our answer is `yes, we intend to introduce ordinals, where $\varphi_0$ is the first one'.

Let $f$, $g$, $a$, and $c$ be pairwise distinct functions where: {\sc (i)}$f(a) = g$, $f(g) = c$, $f(c) = c$, $f(f) = f$, and $f(r) = \underline{\mathfrak{0}}$ for the remaining values $r$; {\sc (ii)} $g(a) = a$, $g(g) = g$, and $g(r) = \underline{\mathfrak{0}}$ for the remaining values $r$; {\sc (iii)} $a$ and $c$ are arbitrary, as long they are neither $\underline{\mathfrak{0}}$ nor $\underline{\mathfrak{1}}$. Then, $(f\circ g)\circ f = \varphi_0$, while $f\circ (g\circ f) = l$, where $l(a) = c$, $l(l) = l$, and $l(r) = \underline{\mathfrak{0}}$ for the remaining values $r$. Thus, $\mathfrak{F}$-composition among functions which act on each other, if they exist, is not necessarily associative.

\begin{teorema}\label{associatividadeentrenaocirculares}
Let $f$, $g$, and $h$ be terms such that neither one of them acts on the remaining ones. Then $\mathfrak{F}$-Composition is associative among $f$, $g$, and $h$.
\end{teorema}

\begin{description}
\item[Proof] Both $f\circ (g \circ h) = p$ and $(f\circ g) \circ h = q$ correspond (Definition \ref{composicao}) to $y = f(g(h(x)))$, if $x$ is different of $f$, $g$, $h$, $p$, and $q$. So, $f\circ (g \circ h) = (f\circ g) \circ h$ for those values of $x$. But $f$, $g$, and $h$ do not act on $f$, $g$, or $h$; and no $\mathfrak{F}$-composition acts on any of its factors. So, $f\circ (g \circ h) = (f\circ g) \circ h = f(g(h(x)))$.
\end{description}

Evaluation $f_1^2$ is not associative as well. Consider, e.g., $x(\underline{\mathfrak{1}}(x))$, for $x$ different of $\underline{\mathfrak{0}}$ and different of $\underline{\mathfrak{1}}$, and such that $x(\underline{\mathfrak{1}}) = \underline{\mathfrak{0}}$ (functions like that do exist, as we can see later on). Thus, $x(\underline{\mathfrak{1}}(x)) = x(x) = x$. If evaluation was associative, we would have $x(\underline{\mathfrak{1}}(x)) = x(\underline{\mathfrak{1}})(x) = \underline{\mathfrak{0}}(x) = \underline{\mathfrak{0}}$; a contradiction, since we assumed $x\neq\underline{\mathfrak{0}}$! Of course this {\em rationale\/} works only if we prove the existence of a function like $x$ and that $\underline{\mathfrak{0}}\neq\underline{\mathfrak{1}}$. Both claims are the subject of the next two theorems.

It is good news that evaluation is not associative. According to {\sc F2}, for all $t$, $g(t) = (g(g))(t)$, since $g(g) = g$. If evaluation was associative, we would have $g(t) = g(g(t))$ and, thus, $g = g\circ g$. So, $\mathfrak{F}$-composition would be idempotent.

\begin{teorema}

$\underline{\mathfrak{0}}\neq\underline{\mathfrak{1}}$.

\end{teorema}

\begin{description}
\item[Proof] $\underline{\mathfrak{0}}\circ\underline{\mathfrak{0}} = \varphi_0\neq\underline{\mathfrak{0}}$ (Theorem \ref{primeirophizero}). But $\underline{\mathfrak{0}}(\varphi_0) = \underline{\mathfrak{0}}$, while $\underline{\mathfrak{1}}(\varphi_0) = \varphi_0$. Thus, from Theorem \ref{igualdadefuncoes}, $\underline{\mathfrak{0}}\neq\underline{\mathfrak{1}}$.
\end{description}

\begin{teorema}\label{segundophizero}

For any $x$ we have $\underline{\mathfrak{0}}\circ x = x\circ\underline{\mathfrak{0}} = \varphi_0$

\end{teorema}

\begin{description}
\item[Proof] First we prove $\underline{\mathfrak{0}}\circ x = \varphi_0$. From Definition \ref{composicao} and axiom {\sc F5}, $(\underline{\mathfrak{0}}\circ x)(t) = \underline{\mathfrak{0}}(x(t)) = \underline{\mathfrak{0}}$ for any $t\neq \underline{\mathfrak{0}}$ and $t\neq x$. But {\sc F5} demands $\underline{\mathfrak{0}}\circ x$ is different of $\underline{\mathfrak{0}}$. Hence, $(\underline{\mathfrak{0}}\circ x)(\underline{\mathfrak{0}}) = \underline{\mathfrak{0}}$. Regarding $x$, there are three possibilities: ({\sc i}) $x = \underline{\mathfrak{0}}$; ({\sc ii}) $x = \varphi_0$; ({\sc iii}) $x$ is neither $\underline{\mathfrak{0}}$ nor $\varphi_0$. The first case corresponds to Theorem \ref{primeirophizero}, which entails $\underline{\mathfrak{0}}\circ x = \varphi_0$. In the second case, if $\underline{\mathfrak{0}}\circ x$ is different of $x = \varphi_0$, then $(\underline{\mathfrak{0}}\circ x)(\varphi_0) = \underline{\mathfrak{0}}$. But that would entail $(\underline{\mathfrak{0}}\circ x)(t) = \underline{\mathfrak{0}}$ for any $t\neq \underline{\mathfrak{0}}\circ x$, which corresponds exactly to function $\varphi_0$ proven in Theorem \ref{primeirophizero}, a contradiction. So, $\underline{\mathfrak{0}}\circ x$ is indeed $\varphi_0$, when $x = \varphi_0$. Concerning last case, since $x\neq\underline{\mathfrak{0}}$ and $x\neq\varphi_0$, then (Theorem \ref{igualdadefuncoes}) there is $t\neq\varphi_0$ such that $x(t)\neq\underline{\mathfrak{0}}$ for $x\neq t$. From {\sc F5}, $(\underline{\mathfrak{0}}\circ x)(t) = \underline{\mathfrak{0}}$ for such value of $t$. But once again we have a function $\underline{\mathfrak{0}}\circ x$ such that $(\underline{\mathfrak{0}}\circ x)(t) = \underline{\mathfrak{0}}$ for any $t\neq \underline{\mathfrak{0}}$, which corresponds to $\varphi_0$ from Theorem \ref{primeirophizero}. Concerning the identity $x\circ\underline{\mathfrak{0}} = \varphi_0$, the proof is analogous. If $h = x\circ\underline{\mathfrak{0}}$, then $h\neq\underline{\mathfrak{0}}$ (from {\sc F5}). Besides, for any $x$, $x\neq h$ entails $h(x) = \underline{\mathfrak{0}}$. Hence, $h = \varphi_0$.
\end{description}

\begin{teorema}\label{fechamentorenato}

$\forall x(x(\underline{\mathfrak{0}}) = \underline{\mathfrak{0}})$.

\end{teorema}

\begin{description}
\item[Proof] $x\circ \underline{\mathfrak{0}} = \varphi_0$ (Theorem \ref{segundophizero}). Hence, for any $t$, $(t\neq x\wedge t\neq\underline{\mathfrak{0}}\wedge t\neq\varphi_0)\Rightarrow \varphi_0(t) = x(\underline{\mathfrak{0}}(t)) = x(\underline{\mathfrak{0}})$. But $\varphi_0(t) = \underline{\mathfrak{0}}$ for any $t\neq\varphi_0$. Thus, $x(\underline{\mathfrak{0}}) = \underline{\mathfrak{0}}$.
\end{description}

\begin{teorema}\label{umbarracompostocomumbarraehumbarra}

$\underline{\mathfrak{1}}\circ\underline{\mathfrak{1}} = \underline{\mathfrak{1}}$.

\end{teorema}

\begin{description}
\item[Proof] From {\sc F5}, $h = \underline{\mathfrak{1}}\circ\underline{\mathfrak{1}}$ entails $h(x) = x$ for any $x\neq\underline{\mathfrak{1}}$. Since the $\mathfrak{F}$-composition is unique and $\underline{\mathfrak{1}}$ guarantees all demanded conditions, then $h = \underline{\mathfrak{1}}$.
\end{description}

\begin{teorema}\label{dezenovedemarco}

$\forall x \forall y(x\circ y = \underline{\mathfrak{1}}\Rightarrow (x = \underline{\mathfrak{1}} \wedge y = \underline{\mathfrak{1}}))$.

\end{teorema}

\begin{description}
\item[Proof] Let $x\neq\underline{\mathfrak{1}}$. Then $\underline{\mathfrak{1}}(x) = \underline{\mathfrak{0}}$, according to {\sc F5}. But that happens only for $x = \underline{\mathfrak{0}}$. And $\underline{\mathfrak{0}}\circ y\neq\underline{\mathfrak{1}}$. Analogous argument holds for $y\neq\underline{\mathfrak{1}}$.
\end{description}

There can be no functions different of $\underline{\mathfrak{1}}$ such that their composition is $\underline{\mathfrak{1}}$.

Next theorem is important for a better understanding about {\sc F1}, although its proof does not demand the use of such a postulate.

\begin{teorema}\label{composicaodexcomumbarra}

$\forall x((x\neq \underline{\mathfrak{0}} \wedge x(\underline{\mathfrak{1}}) = \underline{\mathfrak{0}})\Leftrightarrow (\underline{\mathfrak{1}}\circ x = x\wedge x\circ\underline{\mathfrak{1}} = x \wedge x\neq\underline{\mathfrak{1}}))$.

\end{teorema}

\begin{description}
\item[Proof] The $\Rightarrow$ part. If $x(\underline{\mathfrak{1}}) = \underline{\mathfrak{0}}$, then $x\neq\underline{\mathfrak{1}}$, since $\underline{\mathfrak{1}}(\underline{\mathfrak{1}}) = \underline{\mathfrak{1}}$. If $x\circ \underline{\mathfrak{1}} = h$, then, for any $t$ different of $x$, $\underline{\mathfrak{1}}$, and $h$, we have $h(t) = x(\underline{\mathfrak{1}}(t)) = x(t)$. If $h = x$, then $h$ satisfies all conditions from {\sc F5}, since $x\neq\underline{\mathfrak{1}}$, $h(\underline{\mathfrak{1}}) = x(\underline{\mathfrak{1}}) = \underline{\mathfrak{0}}$, and $h(t) = x(t)$ for any $t\neq h$. Since {\sc F5} demands $h$ to be unique, then $h = x$. If $\underline{\mathfrak{1}}\circ x = h$, we use an analogous argument. For the $\Leftarrow$ part, $\underline{\mathfrak{1}}\circ x = x\circ\underline{\mathfrak{1}} = x$ entails $x\neq\underline{\mathfrak{0}}$ (Theorem \ref{segundophizero}). Since $x\neq\underline{\mathfrak{1}}$, then {\sc F5} demands for the $\mathfrak{F}$-composition $x$ that $x(\underline{\mathfrak{1}}) = \underline{\mathfrak{0}}$.
\end{description}

\begin{observacao}\label{observaunicidade}

If it wasn't for the uniqueness of $\mathfrak{F}$-compositions in {\sc F5}, Flow would be consistent with the existence of many functions, like $x$ and $h$ (from Theorem \ref{composicaodexcomumbarra}), which ``do'' the same thing. We refer to such functions as {\em clones\/}. For a brief investigation about clones see Definition \ref{definesim} and its subsequent discussion. Clones are meant to be different functions $x$ and $h$ which share the same images $x(t)$ and $h(t)$ for any $t$ different of both $x$ and $h$, and such that $x(h) = \underline{\mathfrak{0}}$ and $h(x) = \underline{\mathfrak{0}}$. From Theorem \ref{igualdadefuncoes}, $x\neq h$ (recall {\sc F2}). We use this opportunity to prove Theorem \ref{primeirophizero}, since $\underline{\mathfrak{0}}\circ\underline{\mathfrak{0}} = \varphi_0$, where $\varphi_0$ and $\underline{\mathfrak{0}}$ are clones. But we cease using clones when we talk about other functions. That is why we refer to {\sc F1} as ``weak extensionality''. A strong extensionality postulate would demand that other clones besides $\underline{\mathfrak{0}}$ and $\varphi_0$ cannot exist. We do something like that in some remaining postulates where we use the quantifier $\exists!$.

Another issue concerning $\mathfrak{F}$-composition is its non-associativity. The way we introduced `composition' is not the only one possible. We discuss about that at the end of this paper.
\end{observacao}

\begin{definicao}\label{defineSigmasigma}

Given a term $f$, $g$ is the $\mathfrak{F}$-successor of $f$, and we denote this by $\Sigma(f,g)$, iff $f(g) = \underline{\mathfrak{0}} \wedge \forall x(x\neq g \Rightarrow g(x) = f(x))$.

\end{definicao}

For example, let $f$ by given by $f(a) = b$, $f(b) = c$, $f(c) = c$, $f(f) = f$, and $f(r) = \underline{\mathfrak{0}}$ for the remaining values $r$. Now let $g$ be given by $g(a) = b$, $g(b) = c$, $g(c) = c$, $g(f) = f$, $g(g) = g$, and $g(r) = \underline{\mathfrak{0}}$ for the remaining values $r$, where $g\neq f$. In that case, $\Sigma(f,g)$, if functions like those exist. Observe this has nothing to do with our previous discussion about clones.

\begin{description}

\item[\sc F6 - $\mathfrak{F}$-Successor Function] $\exists!\sigma(\sigma\neq\underline{\mathfrak{0}}\wedge \forall f((f\neq\sigma\wedge f\neq\underline{\mathfrak{0}})\Rightarrow \\(\exists g(\Sigma(f,g)\Leftrightarrow \sigma(f) = g)\vee (\forall h(\neg\Sigma(f,h))\Leftrightarrow \sigma(f) = \underline{\mathfrak{0}}))))$.

\end{description}

This last axiom states the existence and uniqueness of a special function $\sigma$. From now on we write mostly $\sigma_f$ for $\sigma(f)$. Every time we use the symbol $\sigma$ we are referring to the same term from {\sc F6}: the only one such that $\sigma_f = g \neq \underline{\mathfrak{0}}$ is equivalent to $\Sigma(f,g)$, and $\sigma_f = \underline{\mathfrak{0}}$ is equivalent to $\forall h(\neg\Sigma(f,h))$, as long $f$ is neither $\sigma$ nor $\underline{\mathfrak{0}}$. Thus, $\sigma$ successfully `signals' $\mathfrak{F}$-successors, when they exist, for all functions, except when $f$ is either $\underline{\mathfrak{0}}$ or $\sigma$. In those cases, we have $\Sigma(\underline{\mathfrak{0}},\varphi_0)$ and $\sigma_{\underline{\mathfrak{0}}} = \underline{\mathfrak{0}}$ (Theorem \ref{fechamentorenato}), while $\forall h (\neg\Sigma(\sigma, h))$ and $\sigma_{\sigma} = \sigma$ ({\sc F2}).

Next we want to grant the existence of $\mathfrak{F}$-successors for many other terms. As previously announced, this axiom states the existence of an hierarchy of functions, where $\varphi_0$ is the first one.

\begin{description}

\item[\sc F7 - Infinity] $\exists i ((\forall t (i(t) = t \vee i(t) = \underline{\mathfrak{0}}))\wedge \sigma_i \neq \underline{\mathfrak{0}} \wedge (i(\varphi_0) = \varphi_0 \wedge \forall x ((x\neq\underline{\mathfrak{0}}\wedge i(x) = x) \Rightarrow (i(\sigma_x) = \sigma_x \wedge \sigma_x\neq\underline{\mathfrak{0}}))))$.

\end{description}

\begin{definicao}\label{funcaoindutiva}
Any $i$ which satisfies {\sc F7} is said to be {\em inductive\/}.
\end{definicao}

Since the existence of $\varphi_0$ is granted by {\sc F5}, we can use $\sigma$ to justify the existence of $\varphi_1 = \sigma_{\varphi_0}$ such that $\varphi_1(\varphi_1) = \varphi_1$, $\varphi_1(\varphi_0) = \varphi_0$, and for the remaining values $r$ (those who are neither $\varphi_0$ nor $\varphi_1$) we have $\varphi_1(r) = \underline{\mathfrak{0}}$. That happens because {\sc F7} states the existence of at least one other function $i$ and infinitely many other functions. It says $i(\varphi_0) = \varphi_0$. Besides, there is a non-$\underline{\mathfrak{0}}$ $\mathfrak{F}$-successor of $\varphi_0$ such that $i(\sigma_{\varphi_0}) = \sigma_{\varphi_0}$. More than that, if $x$ admits a non-$\underline{\mathfrak{0}}$ $\mathfrak{F}$-successor $\sigma_x$ (where $i(x) = x$), then $i(\sigma_x) = \sigma_x$, where $\sigma_x \neq \underline{\mathfrak{0}}$. Thus, $\varphi_0\neq\underline{\mathfrak{0}}$ and $\sigma_{\varphi_0} = \varphi_1$, where $\varphi_1\neq\varphi_0$ and $\varphi_1\neq\underline{\mathfrak{0}}$. Analogously we can get $\varphi_2$, $\varphi_3$, and so on. Along with those terms $\varphi_n$, {\sc F7} says any inductive function $i$ admits its own non-$\underline{\mathfrak{0}}$ $\mathfrak{F}$-successor $\sigma_i$.

Subscripts $0$, $1$, $2$, $3$, etc., are metalinguistic symbols based on an alphabet of ten symbols (the usual decimal numeral system) which follows the lexicographic order $\prec$, where $0\prec 1\prec 2\cdots\prec 8\prec 9$. If $n$ is a subscript, then $n+1$ corresponds to the next subscript, in accordance to the lexicographic order. In that case, we write $n\prec n+1$. $n+m$ is an abbreviation for $(...(...((n+1)+1)+...1)...)$ with $m$ occurrences of $+$ and $m$ occurrences of pairs of parentheses. Again we have $n\prec n+m$. Besides, $\prec$ is a strict total order. That fact allows us to talk about a minimum value between subscripts $m$ and $n$: $\mbox{min}\{ m, n\}$ is $m$ iff $m\prec n$, it is $n$ iff $n\prec m$, and it is either one of them if $m = n$. Of course, $m = n$ iff $\neg(m\prec n) \wedge \neg (n\prec m)$. If $m\prec n \vee m = n$, we denote this by $m\preceq n$. Such a vocabulary of ten symbols endowed with $\prec$ is called here {\em (meta) language\/} $\mathcal L$.

{\sc F7} provides us some sort of ``recursive definition'' for functions $\varphi_n$, while it allows as well to guarantee the existence of inductive functions: {\sc (i)} $\varphi_0$ is such that $\varphi_0(x)$ is $\varphi_0$ if $x = \varphi_0$ and $\underline{\mathfrak{0}}$ otherwise; {\sc (ii)} $\varphi_{n+1}$ is such that $\varphi_{n+1}(\varphi_{n+1}) = \varphi_{n+1}$, $\varphi_{n+1}\neq \varphi_n$, and $\varphi_{n+1}(x) = \varphi_n(x)$ for any $x$ different of $\varphi_{n+1}$.

Observe that $\varphi_{n+1}(\varphi_n) = \varphi_n(\varphi_n) = \varphi_n$, while $\varphi_n(\varphi_{n+1}) = \underline{\mathfrak{0}}$. Moreover, $\varphi_{n+2}(\varphi_{n+1}) = \varphi_{n+1}$, and $\varphi_{n+2}(\varphi_n) = \varphi_{n+1}(\varphi_n) = \varphi_n$; while $\varphi_n(\varphi_{n+2}) = \underline{\mathfrak{0}}$.

{\sc Figure 1} illustrates how to represent some functions $f$ in a quite intuitive way. A {\em diagram\/} of $f$ is formed by a rectangle. On the left top corner inside the rectangle we find label $f$. The remaining labels refer to terms $x$ such that either $f(x)\neq \underline{\mathfrak{0}}$ or there is $t$ such that $f(t) = x$. For each $x$ inside the rectangle there is a unique corresponding arrow which indicates the image of $x$ by $f$, as long $x$ is not $f$ itself. Due to self-reference, label $f$ at the left top corner inside the rectangle does not need to be attached to any arrow.

From left to right, the first diagram refers to $\varphi_0$. It says, for any $x$, $\varphi_0(x)$ is $\underline{\mathfrak{0}}$, except for $\varphi_0$ itself. The second diagram says $\varphi_1(\varphi_1) = \varphi_1$, and $\varphi_1(\varphi_0) = \varphi_0$. The circular arrow associated to $\varphi_0$ in the second diagram says $\varphi_1(\varphi_0) = \varphi_0$. Self-reference postulate does not need to be represented in diagrams. The third diagram says $\varphi_2(\varphi_2) = \varphi_2$, $\varphi_2(\varphi_1) = \varphi_1$, and $\varphi_2(\varphi_0) = \varphi_0$. Finally, for the sake of illustration, the last diagram corresponds to an $f$ such that $f(a) = b$, $f(b) = c$, and $f(c) = \underline{\mathfrak{0}}$. That is why there is no arrow `starting' at $c$. In those cases we represent $c$ outside the rectangle. The existence of functions like that is granted by {\sc F10}$_{\alpha}$. The diagrams of $\underline{\mathfrak{0}}$ and $\underline{\mathfrak{1}}$ are, respectively, a blank rectangle and a filled in black rectangle. Other examples are provided in the next paragraphs.\\

\begin{picture}(60,50)

\put(3,13){\framebox(50,40)}


\put(10,43){$\varphi_0$}


\put(63,13){\framebox(50,40)}


\put(70,43){$\varphi_1$}

\put(70,26){$\boldmath \curvearrowleft$}

\put(70,20){$\varphi_0$}


\put(123,13){\framebox(50,40)}


\put(130,43){$\varphi_2$}

\put(130,26){$\boldmath \curvearrowleft$}

\put(130,20){$\varphi_0$}

\put(160,26){$\boldmath \curvearrowleft$}

\put(160,20){$\varphi_1$}

\put(190,30){$\cdots$}

\put(250,13){\framebox(50,40)}

\put(255,43){$f$}

\put(255,23){$a$}

\put(275,23){$b$}

\put(315,23){$c$}

\put(282,26){\vector(2,0){30}}

\put(262,26){\vector(2,0){10}}

\end{picture}

\vspace{-3mm}

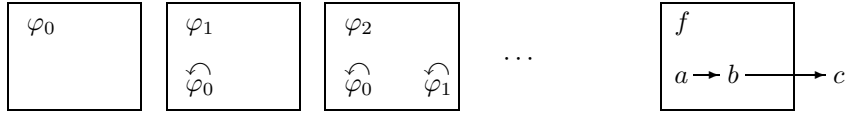
{\sc Figure 1:} From left to right, diagrams of $\varphi_0$, $\varphi_1$, $\varphi_2$, and an arbitrary $f$.

Observe that $\varphi_{m+n}(\varphi_n) = \varphi_n$, for any $m$ and $n$ of $\mathcal L$.

\begin{definicao}\label{actionaction}

$f[t]$ iff $t\neq f \wedge f(t) \neq \underline{\mathfrak{0}}$. We read $f[t]$ as ``{\em $f$ acts on $t$\/}''.

\end{definicao}

Whereas $f(t)$ is a term, $f[t]$ abbreviates a formula. No $f$ acts on itself. The intuitive idea is to allow us to talk about what a function $f$ {\em does\/}. For example, $\underline{\mathfrak{0}}$ and $\varphi_0$ do nothing whatsoever, since there is no $t$ on which they act. On the other hand, there is a term $t$ on which $\varphi_1$ acts, namely, $\varphi_0$.

\begin{teorema}\label{noaction}

$\underline{\mathfrak{0}}$ and $\varphi_0$ are the only functions who do not act on any $t$.

\end{teorema}

\begin{description}
\item[Proof] If $f$ does not act on any $t$, then $\forall t (t = f\vee f(t) = \underline{\mathfrak{0}})$. Suppose $f$ is neither $\underline{\mathfrak{0}}$ nor $\varphi_0$. Then there is $t$ such that $t\neq f \wedge f(t)\neq\underline{\mathfrak{0}}$. But that entails $f[t]$. So, the only functions which do not act on any $t$ are $\underline{\mathfrak{0}}$ and $\varphi_0$.
\end{description}

\begin{teorema}\label{sucessordeumbarra}
$\sigma_{\underline{\mathfrak{1}}} = \underline{\mathfrak{0}}$.
\end{teorema}

\begin{description}
\item[Proof] $\sigma_{\underline{\mathfrak{1}}}\neq\underline{\mathfrak{0}}\Rightarrow\underline{\mathfrak{1}}(\sigma_{\underline{\mathfrak{1}}}) = \underline{\mathfrak{0}}$ (Definition \ref{defineSigmasigma}). That happens only if $\sigma_{\underline{\mathfrak{1}}} = \underline{\mathfrak{0}}$.
\end{description}

\begin{teorema}\label{psi}

For any $g$, $\sigma_g \neq \underline{\mathfrak{1}}$

\end{teorema}

\begin{description}
\item[Proof] Suppose there is $g$ such that $\sigma_g = \underline{\mathfrak{1}}$. Since $\underline{\mathfrak{1}}\neq\underline{\mathfrak{0}}$, then $g(\underline{\mathfrak{1}}) = \underline{\mathfrak{0}}$ and $\underline{\mathfrak{1}}(g) = g$, according to Definition \ref{defineSigmasigma}. So, $g\neq\underline{\mathfrak{1}}$. Once again from Definition \ref{defineSigmasigma}, $g(x) = \underline{\mathfrak{1}}(x)$ for any $x\neq\underline{\mathfrak{1}}$. But those are the same conditions for $\underline{\mathfrak{1}}\circ\underline{\mathfrak{1}} = g$, according to {\sc F5}. Since any $\mathfrak{F}$-composition is unique, then $g = \underline{\mathfrak{1}}$ (Theorem \ref{umbarracompostocomumbarraehumbarra}), which contradicts the assumption $g\neq\underline{\mathfrak{1}}$. So, there is no such $g$.
\end{description}

In other words, the existence of some `functions' in Flow is forbidden.

\begin{definicao}\label{restricaodefuncao}
For any function $f$, a {\em restriction\/} $g$ {\em of\/} $f$ is defined as

\noindent
$$g\subseteq f \;\;\mbox{iff}\;\; g\neq\underline{\mathfrak{0}}\wedge\forall x((g[x]\Rightarrow f[x]) \wedge ((g[x]\wedge f[x])\Rightarrow f(x) = g(x))),$$

\end{definicao}

Proper restrictions are defined as $g\subset f\;\;\mbox{iff}\;\; g\subseteq f \wedge g\neq f$. We abbreviate $\neg (g\subseteq f)$ and $\neg (g\subset f)$ as, respectively, $g\not\subseteq f$ and $g\not\subset f$. For example, $\varphi_1\subseteq \varphi_3$, $\varphi_1\subset \varphi_3$, and $\varphi_3\not\subseteq \varphi_1$.

\begin{teorema}\label{restricaoderenato}

$\forall f \forall g (g\subset f \Rightarrow g(f) = \underline{\mathfrak{0}})$.

\end{teorema}

\begin{description}
\item[Proof] Suppose $g(f)\neq\underline{\mathfrak{0}}$. Since $g\neq f$, then $g[f]$. But for any $f$ we have $\neg f[f]$. In other words, $\neg (g[f]\Rightarrow f[f])$. Hence, $g\not\subseteq f$.
\end{description}

Some of the most useful restrictions are obtained from a given formula $F$, in a way which resembles the well known Separation Scheme. That is achieved thanks to careful considerations regarding the $\mathfrak{F}$-successor function $\sigma$. But before that, we need more concepts, since we are interested on a vast number of situations.

\begin{definicao}\label{abrangente}

A term $f$ is {\em comprehensive\/} iff there is $g$ such that $g\neq\underline{\mathfrak{0}}$, $g\subseteq f$, and $\sigma_g = \underline{\mathfrak{0}}$; and we denote that by $\mathbb{C}(f)$. Otherwise, $f$ is {\em uncomprehensive\/}.

\end{definicao}

Comprehensive functions are supposed to describe ``huge'' functions who act on ``many terms'', and they are partially regulated by some of the the next axioms. If $f$ itself is such that $\sigma_f = \underline{\mathfrak{0}}$, then $f$ is comprehensive, except when $f = \underline{\mathfrak{0}}$. Besides, it is easy to see that $\underline{\mathfrak{1}}$ is comprehensive and $\underline{\mathfrak{0}}$ is uncomprehensive. On the other hand, there are other comprehensive functions which, by the way, play an important role within our proposal for cardinalities and model theory. But for now we are mostly interested on uncomprehensive functions, as it follows.

\begin{definicao}\label{definindoemergente}

$\mathbb{E}(f)$ iff {\sc (i)} $\sigma_f\neq\underline{\mathfrak{0}}$; {\sc (ii)} $\forall x(f[x]\Rightarrow (\mathbb{E}(x)\wedge \mathbb{E}(f(x))))$.

\end{definicao}

We read $\mathbb{E}(f)$ as ``$f$ is {\em emergent\/}''. $\varphi_0$ is vacuously emergent. That entails $\varphi_1$ is emergent. A function like $f$ given by $f(\varphi_0) = \varphi_1$, $f(\varphi_1) = \varphi_0$, $f(f) = f$, and $f(r) = \underline{\mathfrak{0}}$ (where $r$ stands for the remaining values $r$) is emergent. Obviously, last claim demands the existence of such an $f$ and the existence of $\sigma_f\neq\underline{\mathfrak{0}}$, something which is accomplished thanks to the next postulates. Emergent functions are supposed to be uncomprehensive terms who act only on uncomprehensive terms.

\begin{definicao}\label{lurking}

$g\trianglelefteq f$ iff $g\neq\underline{\mathfrak{0}}\wedge\forall x(g[x]\Rightarrow ((f[x]\vee \exists a(f(a) = x))\wedge (f[g(x)]\vee \exists b (f(b) = g(x)))))$. If $g\trianglelefteq f$, we say $g$ {\em lurks\/} $f$. Besides, $g\triangleleft f$ iff $g\trianglelefteq f$ and $g\neq f$. In that case we say $g$ {\em properly lurks\/} $f$. Finally, $\neg (g\trianglelefteq f)$ and $\neg (g\triangleleft f)$ are abbreviated as $g\ntrianglelefteq f$ and $g\ntriangleleft f$, respectively.

\end{definicao}

As an example, if $f$ is the same function of last paragraph, then $f\triangleleft \varphi_2$ and $\varphi_2\triangleleft f$. As another example, $\varphi_2\triangleleft \varphi_4$, but $\varphi_4\ntrianglelefteq \varphi_2$.

\begin{teorema}

If $g\subseteq f$ and $g\neq\underline{\mathfrak{0}}$, then $g\trianglelefteq f$.

\end{teorema}

The proof is straightforward. The converse is obviously not valid.

\begin{definicao}\label{definindomaximapotencia}

$h = \mathfrak{p}(f)$ iff $\forall x (h[x]\Leftrightarrow x\trianglelefteq f)$. $h$ is the {\em full power\/} of $f$.

\end{definicao}

The full power $\mathfrak{p}(f)$ of $f$ acts on all terms $x$ who lurk $f$.

\begin{teorema}

If $f\subset\underline{\mathfrak{1}}$ acts on $n$ terms, then $\mathfrak{p}(f)$ acts on $(n+1)^n$ terms.

\end{teorema}

\begin{description}
\item[Proof] If $f$ acts on $n$ terms $x_i$ and $g$ lurks $f$, each $x_i$ may correspond to any $x_j$ (where eventually $x_j = x_i$) in the sense we may have $g(x_i) = x_j$. On the other hand, we may have $g(x_i) = \underline{\mathfrak{0}}$ as well. Thus, for each $x_i$ ($n$ possible values) there are $n+1$ possible images. So, there are $(n+1)^n$ terms who lurk $f$.
\end{description}

\begin{description}

\item[\sc F8 - Coherence] $\forall f(\mathbb{E}(f)\Rightarrow (\forall g ((g\trianglelefteq f\Rightarrow \sigma_g\neq\underline{\mathfrak{0}})\wedge \forall h ((h[g]\Leftrightarrow g\trianglelefteq f)\Rightarrow \sigma_h\neq\underline{\mathfrak{0}}))\wedge \forall i(\forall t(i[t]\Rightarrow \exists x(x[t]\wedge f[x])))\Rightarrow \sigma_i\neq\underline{\mathfrak{0}}))$.

\end{description}

If $f$ is emergent, then any $g$ (if it exists) who lurks $f$ has a non-$\underline{\mathfrak{0}}$ $\mathfrak{F}$-successor. Function $h$ at {\sc F8} refers to the full power of $f$, if it exists. Finally, term $i$ is useful for dealing with arbitrary union, if we are able to define it. The usual power of a ZF-set (Definition \ref{definindowp} and Theorem \ref{potenciaZF}) and arbitrary unions (Definition \ref{definindouniaoarbitraria} and Theorem \ref{uniaoarbitrariadezfsets}) are definable thanks to the next postulate coupled with {\sc F8}.

Observe {\sc F8} is an existence postulate. For understanding this, recall {\sc F6}, which provides necessary and sufficient conditions for knowing if $\sigma_f\neq\underline{\mathfrak{0}}$: {\em there must be\/} a $g$ different of $f$ such that certain conditions are met. The point here is that {\sc F8} grants the existence of certain terms which are $\mathfrak{F}$-successors of others, as long some conditions are met. When we say, as above, that $\sigma_g\neq\underline{\mathfrak{0}}$, we state {\em there is\/} a function $z\neq\underline{\mathfrak{0}}$ such that $\sigma_g = z$. The same happens to $\sigma_h$ and $\sigma_i$.

\begin{definicao}\label{defineZFconjunto}

$\mathbb{Z}(f)$ iff {\sc (i)} $\sigma_f\neq\underline{\mathfrak{0}}$; {\sc (ii)} $\forall x (f[x]\Rightarrow (f(x) = x\wedge \mathbb{Z}(x)))$.

\end{definicao}

$\mathbb{Z}(f)$ reads ``$f$ is a ZF-set''. Any ZF-set $f$ is a restriction of $\underline{\mathfrak{1}}$, where $\sigma_f\neq\underline{\mathfrak{0}}$; and if $f$ acts on $x$, then $\sigma_x \neq \underline{\mathfrak{0}}$. Besides, every restriction of a ZF-set has its own non-$\underline{\mathfrak{0}}$ $\mathfrak{F}$-successor ({\sc F8}). So, every restriction of a ZF-set is a ZF-set. Observe $\varphi_0$ vacuously satisfies last definition. Therefore, $\varphi_1$ is a ZF-set and so on. Observe as well, out of curiosity, $\sigma$ is not a ZF-set, since $\sigma(\varphi_0) = \varphi_1 \neq \varphi_0$. Every ZF-set is an emergent function which is a restriction of $\underline{\mathfrak{1}}$.

\begin{teorema}

No ZF-set is comprehensive.

\end{teorema}

\begin{description}
\item[Proof] If $f$ is a ZF-set, then $\sigma_f\neq\underline{\mathfrak{0}}$ and $f$ acts on terms $t$ such $\sigma_t\neq\underline{\mathfrak{0}}$. From {\sc F8}, any restriction $g$ of $f$ is such that $\sigma_g\neq\underline{\mathfrak{0}}$.
\end{description}

\begin{teorema}\label{todophinehZFconjunto}

Every $\varphi_n$ is a ZF-set.

\end{teorema}

\begin{description}
\item[Proof] $\mathbb{Z}(\varphi_0)$ is vacuously valid. Now, let $n>0$. Then any $\varphi_n$ acts only on $\varphi_m$ and $\varphi_n(\varphi_m) = \varphi_m$, where $0<m<n$. And each $\varphi_m$ has a non-$\underline{\mathfrak{0}}$ $\mathfrak{F}$-successor, from the definition itself for $\varphi_m$. From {\sc F8}, that entails that any restriction $g$ of $\varphi_n$ has a non-$\underline{\mathfrak{0}}$ $\mathfrak{F}$-successor. So, $\mathbb{Z}(\varphi_n)$ for any $n$ from language $\mathcal L$.
\end{description}

\begin{definicao}\label{delimitadorarestricao}

Let $f$ be a function and $F(t)$ be a formula where all occurrences of $t$ are free. We say $g$ is {\em restriction of $f$ under\/} $F(t)$, and denote this by $g = f\big|_{F(t)}$, iff: {\sc (i)} $g\neq\underline{\mathfrak{0}}$; {\sc (ii)} $f(\sigma_g) = \underline{\mathfrak{0}} \vee \neg F(g)$; {\sc (iii)} $g\neq f\Rightarrow g(f) = \underline{\mathfrak{0}}$; {\sc (iv)} $\forall t((t\neq f\wedge t\neq g)\Rightarrow ((g(t) = f(t)\wedge F(t)\wedge f[t])\vee (g(t) = \underline{\mathfrak{0}}\wedge (\neg F(t)\vee \neg f[t]))))$.

\end{definicao}

Formula $g = f\big|_{F(t)}$ is somehow equivalent to $g\subseteq f$, as we see next.

\begin{teorema}

If $g\subseteq f$, then it is possible to state a formula $F$ where $g = f\big|_F$.

\end{teorema}

\begin{description}
\item[Proof] If $g\subseteq f$, assume as formula $F(t)$ the next one: $g[t]$. Item ({\sc i}) of Definition \ref{delimitadorarestricao} is a consequence from Definition \ref{restricaodefuncao}. Item ({\sc ii}) of the same definition is granted thanks to the fact that $\neg g[g]$ ($\neg F(g)$) for any $g$. Item ({\sc iii}) is due to Theorem \ref{restricaoderenato}. Finally, item ({\sc iv}) is granted from Definition \ref{restricaodefuncao}.
\end{description}

\begin{teorema}

If $g = f\big|_F$, then $g\subseteq f$.

\end{teorema}

\begin{description}
\item[Proof] If $g = f\big|_F$, then either $g = \varphi_0$ or $g\neq\varphi_0$. In the first case, the proof is immediate by vacuity ($\varphi_0$ does not act on any term). If $g\neq\varphi_0$, then item ({\sc iv}) of Definition \ref{delimitadorarestricao} demands, for any $t$, $g[t]\Rightarrow (f[t]\wedge g(t) = f(t))$. Thus, $g\subseteq f$.
\end{description}

A natural way of getting some restrictions $g$ of $f$ is through {\sc F9}$_{\mathbb{E}}$ below. Observe as well item ({\sc iv}) of Definition \ref{delimitadorarestricao} takes into account the self-reference postulate, since we demand $t\neq f \wedge t\neq g$. Now, if $F(t)$ is a formula (abbreviated by $F$) where all occurrences of $t$ are free, then the following is an axiom.

\begin{description}

\item[\sc F9$_{\mathbb{E}}$ - $\mathbb{E}$-Restriction] $\forall f((\forall x(F(x)\Rightarrow \mathbb{E}(x))\Rightarrow \exists! g(g = f\big|_F)))$.

\end{description}

Subscript $_{\mathbb{E}}$ highlights a strong commitment to emergent functions. Many restrictions due to this last postulate grant the existence of emergent functions and ZF-sets as well. Next theorem, for example, shows we do not need $\mathfrak{F}$-composition to prove there is $\varphi_0$.

\begin{teorema}\label{teoremasobrerestricaodezerobarra}

$\underline{\mathfrak{0}}\big|_{F(x)} = \varphi_0\big|_{F(x)} = \varphi_0$ if $F(x)\Rightarrow \mathbb{E}(x)$.

\end{teorema}

\begin{description}
\item[Proof] Immediate, since neither $\underline{\mathfrak{0}}$ nor $\varphi_0$ act on any term and no restriction can be $\underline{\mathfrak{0}}$.
\end{description}

\begin{teorema}\label{trioderestricoes}

For any emergent $f$ we have: {\sc (i)} $f\big|_{x\neq x\wedge\mathbb{E}(x)} = \varphi_0$;\\ {\sc (ii)} $f\big|_{x = x \wedge \mathbb{E}(x)} = f$; {\sc (iii)} $f\big|_{x \neq f\wedge \mathbb{E}(x)} = f$; {\sc (iv)} $f\big|_{x = f} = \varphi_0$.

\end{teorema}

\begin{description}
\item[Proof] Item {\sc (i)} is proven by vacuity. {\sc (ii)} takes into account all terms where $f$ acts. About {\sc (iii)} and {\sc (iv)}, recall $f$ plays no role into the calculation of its restriction. All that matters are the terms where $f$ acts.
\end{description}

From {\sc F9}$_{\mathbb{E}}$, there are {\em four\/} possible restrictions $g$ of $\varphi_2$. If $F(x)$ is, for example, ``$x = \varphi_0$'', then $g = \varphi_1$. Accordingly, assume $f = \varphi_2$ in {\sc F9}$_{\mathbb{E}}$. So, consider, e.g., $x = \varphi_0$. Such a value for $x$ is different of $\varphi_2$. Besides, $F(\varphi_0)$. That implies $g(\varphi_0) = \varphi_2(\varphi_0) = \varphi_0$. For all remaining values $x\neq\varphi_0$, we know $g$ does not act on $x$. That means $g$ acts solely on $\varphi_0$. And according to {\sc F7}, that function is supposed to be $\varphi_1$. Observe $\varphi_1$ is allowed to have a free occurrence in $F(x)$, according to our Restriction Axiom. Nevertheless $g$ does not act on $\varphi_1$ in our first example. That means either $g(\varphi_1) = \underline{\mathfrak{0}}$ or $g = \varphi_1$. In this case, we have $g = \varphi_1$. Latter on we define a membership relationship $\in$ (Definition \ref{pertencercomoagir}) where $x\in f$ iff $f[x]$ and some conditions are imposed over $f$ and $x$. That entails we can guarantee that in a translation of ZF's Separation Scheme into Flow's language, any free occurrence of $g$ in $F(x)$ will have no impact (in a precise sense). After all, in this first example $F(x)$ is $x = \varphi_0$, while $g$ is $\varphi_1$. For details see Section \ref{standard}.

Resuming, if $F(x)$ is the formula ``$x = \varphi_0 \vee x = \varphi_1$'', then $g = \varphi_2$. If $F(x)$ is ``$x = x$'', then again $g = \varphi_2$. If $F(x)$ is ``$x \neq x$'', then $g = \varphi_0$. The novelty here, however, happens with the formula $F(x)$ given by ``$x = \varphi_1$''. In that case we have a proper restriction $\gamma$ such that $\gamma\neq \varphi_1$, $\gamma(\gamma) = \gamma$, $\gamma(\varphi_1) = \varphi_1$, and $\gamma(x) = \underline{\mathfrak{0}}$ for any $x$ different of $\varphi_1$ and $\gamma$ itself. Thus, $\gamma$ is a {\em new function\/} whose existence is granted thanks to {\sc F9}$_{\mathbb{E}}$ and no other previous postulate.

\begin{definicao}\label{definindowp}

$z$ is the {\em restricted power of\/} $f$ iff $\forall x (x\neq z \Rightarrow((z(x) = x \Leftrightarrow x \subseteq f)\wedge (z(x) = \underline{\mathfrak{0}} \Leftrightarrow x\not\subseteq f)))$. We denote $z$ as $\wp(f)$.

\end{definicao}

For example, $\wp(\varphi_0) = \varphi_1$, $\wp(\varphi_1) = \varphi_2$, and $\wp(\varphi_2) = f$, where $f(\varphi_0) = \varphi_0$, $f(\varphi_1) = \varphi_1$, $f(\varphi_2) = \varphi_2$, $f(\gamma) = \gamma$, $f(f) = f$, and $f(x) = \underline{\mathfrak{0}}$ for the remaining values $x$. Recall $\gamma$ acts only on $\varphi_1$ and $\gamma(\varphi_1) = \varphi_1$.

Observe $z = \wp(f)$ is a restriction of $\underline{\mathfrak{1}}$, even if $f$ is not. Besides, $\wp$ is {\em not\/} a function, but a metalinguistic symbol which helps us to abbreviate the formula $z = \wp(f)$ given by the definition above. Observe as well $\wp(f)\subseteq \mathfrak{p}(f)$, for any $f$. While $\mathfrak{p}(f)$ refers to a function who acts on all terms that lurk $f$, $\wp(f)$ acts on all terms that lurk $f$ as long they are restrictions of $f$.

If $g = \mathfrak{p}(f)$, we can get $\wp(f)$ in the case where $f$ is a ZF-set, in the sense $\wp(f)$ is the term which acts on all restrictions of $f$. For details see Theorem \ref{potenciaZF}.

\begin{observacao}\label{consideracoessobrerestricao}

In a sense, {\sc F9}$_{\mathbb{E}}$ is similar to the Separation Scheme in ZFC, since it states the existence of a unique $g$ obtained from a given $f$ and a formula $F(x)$. Nevertheless, the role of Separation Scheme in ZFC is not limited to grant the existence of subsets. Thanks to that postulate, ZFC avoids antinomies like Russell's paradox. In our case those antinomies are avoided by means of the simple use of Self-Reference (Observation \ref{Russell}). That is one of the reasons why we do not prohibit free occurrences of $g$ in $F(x)$ (like what happens in ZFC). Actually, if we demanded no free occurrences of $g$ in $F(x)$, we would be unable to obtain some useful restrictions, as we can see in the examples below. Nevertheless, we demand $\forall x ((x\neq g \wedge x\neq f) \Rightarrow ((f(x) = g(x)\wedge F(x))\vee (g(x) = \underline{\mathfrak{0}}\wedge\neg F(x)))$ (Definition \ref{delimitadorarestricao}), which is a weaker condition than the prohibition of occurrences of $g$ in $F(x)$. If we recall Observation \ref{Russell}, we can easily see that, for any $y$ and any $r$, neither $y\big|_{x(x)\neq y}$ nor $y\big|_{x(x)\neq r}$ allow us to get any contradiction in the style of Russell's paradox. Finally, if anyone tries to ``define'' a function $f$ from a formula $F(x)$ without using {\sc F9}$_{\mathbb{E}}$, it is perfectly possible to get a contradiction. For example, we can ``define'' a function $f$ as it follows: $\forall x(f(x) = \varphi_0)$. Since no $f$ acts on $\underline{\mathfrak{0}}$, we have a contradiction, namely, $f(\underline{\mathfrak{0}}) = \varphi_0$. Nevertheless, that would not be a definition at all, but simply a new postulate which is inconsistent with our axioms.

\end{observacao}

\begin{definicao}\label{definesim}

$f\sim g$ iff $\forall t ((f[t]\Leftrightarrow g[t])\wedge ((t\neq f \wedge t \neq g)\Rightarrow f(t) = g(t)))$.

\end{definicao}

If $f \sim g$ but $f\neq g$, we say $f$ and $g$ are {\em clones\/}, as discussed in Observation \ref{observaunicidade}. Regarding their images, clones $f$ and $g$ differ solely on $f$ and $g$: $f(g) = \underline{\mathfrak{0}}$, while $g(g) = g$, and $g(f) = \underline{\mathfrak{0}}$, while $f(f) = f$. It is easy to see that $\sim$ is reflexive.

\begin{teorema}\label{marcionoslibertou}

Let $f$ be emergent. Then, for any $h\neq\underline{\mathfrak{0}}$, $h \sim f$ entails $h = f$.

\end{teorema}

\begin{description}
\item[Proof] From Theorem \ref{trioderestricoes}, $f = f\big|_{t = t\wedge\mathbb{E}(t)}$. Since $h \sim f$, $h = f\big|_{t = t\wedge\mathbb{E}(t)}$ as well. But {\sc F9}$_{\mathbb{E}}$ says any restriction is unique. Thus, $h = f$.
\end{description}

This last theorem grants $\underline{\mathfrak{0}}$ and $\varphi_0$ are the only clones in Flow, at least among emergent functions. That is a useful result for discussing subjects like the `quantity' of restrictions of any $f$.

\begin{teorema}
$\forall f (\varphi_0\subseteq f).$
\end{teorema}

\begin{description}
\item[Proof] Theorem \ref{noaction} states $\varphi_0$ is the only function different of $\underline{\mathfrak{0}}$ who do not act on any term. Thus, Definition \ref{restricaodefuncao} grants formula above is proven by vacuity.
\end{description}

\begin{teorema}

For any $x$ and $y$, $x\subset y$ entails $y\not\subset x$.

\end{teorema}

\begin{description}
\item[Proof] If $x$ and $y$ are, respectively, $\varphi_0$ and $\underline{\mathfrak{0}}$, the proof is straightforward, since $\underline{\mathfrak{0}}$ is never any restriction, according to {\sc F9}$_{\mathbb{E}}$. For the remaining cases, $x\subset y$ entails $x\neq y$ and, for any $t$, $x[t]\Rightarrow (y[t]\wedge x(t) = y(t))$. So, there is $t'$ where $y[t']\wedge\neg x[t']$ or $y[t']\wedge x[t']\wedge x(t')\neq y(t')$. In anyone of those cases we have $y\not\subset x$.
\end{description}

\begin{teorema}\label{onlyyou}

$\sigma$ and $\underline{\mathfrak{0}}$ are the only functions $f$ such that $\sigma_f = f$.

\end{teorema}

\begin{description}
\item[Proof] We already know $\sigma_\sigma = \sigma$ ({\sc F2}) and $\sigma_{\underline{\mathfrak{0}}} = \underline{\mathfrak{0}}$ (Theorem \ref{fechamentorenato}). If $f\neq\underline{\mathfrak{0}}$ and $\sigma_f = \underline{\mathfrak{0}}$, then $\sigma_f\neq f$. If $f\neq\sigma$, $f\neq\underline{\mathfrak{0}}$ and $\sigma_f\neq\underline{\mathfrak{0}}$, then {\sc F6} and Definition \ref{defineSigmasigma} demand $\sigma_f\neq f$.
\end{description}

\begin{teorema}

$\forall f(\mathbb{E}(f)\Rightarrow f = \sigma_f\big|_{x\neq f\wedge\mathbb{E}(x)})$.

\end{teorema}

\begin{description}
\item[Proof] $\mathbb{E}(f)$ entails $\sigma_f\neq\underline{\mathfrak{0}}$. So, $\sigma_f$ acts on all terms where $f$ acts. But $\sigma_f$ acts on just one more term: $f$ itself. Since formula $F(x)$ in the restriction above is ``$x\neq f\wedge\mathbb{E}(x)$'' (observe we follow all demands for $F(x)$ in {\sc F9}$_{\mathbb{E}}$), then $\sigma_f\big|_{x\neq f\wedge\mathbb{E}(x)}$ acts exactly on all terms where $f$ acts. Besides, they share all their images, according to {\sc F9}$_{\mathbb{E}}$. So, $\sigma_f\big|_{x\neq f\wedge\mathbb{E}(x)} = f$, since $f[x]\Rightarrow x\neq f$.
\end{description}

\begin{definicao}\label{orderedpair}
$f$ is an {\em ordered pair\/} $(a,b)$, with both values $a$ and $b$ different of $\underline{\mathfrak{0}}$, iff there are $\alpha$ and $\beta$ such that $\alpha\neq f$, $\beta\neq f$, $\alpha\neq a$, $\beta\neq b$ and
$$f(x) = \left\{ \begin{array}{cl}
\alpha & \mbox{if}\; x = \alpha\\
\beta & \mbox{if}\; x = \beta\\
\underline{\mathfrak{0}} & \mbox{if}\; x\neq f \wedge x\neq\alpha \wedge x\neq \beta
\end{array}
\right.$$
\noindent
where $\alpha(a) = a$, $\alpha(x) = \underline{\mathfrak{0}}$ if $x$ is neither $a$ nor $\alpha$, $\beta(a) = a$, $\beta(b) = b$, $\beta(x) = \underline{\mathfrak{0}}$ if $x$ is neither $a$ nor $b$ or $\beta$.

\end{definicao}

We have two kinds of ordered pairs: those where $\alpha\neq b$ (first kind) and those where $\alpha = b$ (second kind). The diagram of the first kind is as follows:

\begin{picture}(60,65)

\put(3,25){\framebox(100,30)}


\put(10,45){$f$}

\put(10,11){\framebox(40,30)}

\put(40,43){$\boldmath \curvearrowleft$}

\put(15,32){$\alpha$}

\put(55,11){\framebox(40,30)}

\put(85,43){$\boldmath \curvearrowleft$}

\put(60,31){$\beta$}

\put(60,12){$a$}

\put(58,16){$\boldmath \curvearrowleft$}

\put(13,16){$\boldmath \curvearrowleft$}

\put(77,17){$\boldmath \curvearrowleft$}

\put(80,12){$b$}

\put(15,12){$a$}

\put(203,23){\framebox(100,40)}


\put(210,53){$g$}

\put(220,5){\framebox(75,45)}

\put(240,8){\framebox(50,30)}

\put(280,40){$\boldmath \curvearrowleft$}

\put(285,51){$\boldmath \curvearrowleft$}

\put(225,39){$\beta$}

\put(244,28){$\alpha = b$}

\put(280,10){$a$}

\put(278,15){$\boldmath \curvearrowleft$}

\end{picture}

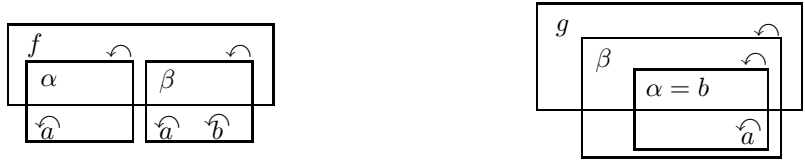
{\sc Figure 2:} Diagrams of $f = (a,b)$ and $g = (a,b)$ of the first and second kind.

Left diagram above concerns the case where $f$ acts only on $\alpha$ and $\beta$, while $\alpha$ acts only on $a$, and $\beta$ acts only on $a$ and $b$. In the particular case where $a = b$, we have $\alpha = \beta$, and the ordered pair $f$ is $(a,a)$. So, $(a,a)$ is a function $f$ which acts solely on $\alpha$, while $\alpha$ acts solely on $a$. Observe that $f(a) = f(b) = \underline{\mathfrak{0}}$ ($f$ never acts neither on $a$ nor on $b$) if $f$ is an ordered pair of the first kind. Thus, $f$ is $(a,b)$ iff $f$ acts only on $\alpha$ and $\beta$, which act, respectively, only on $a$ and only on $a$ and $b$. To get $(b,a)$, all we have to do is to exchange $\alpha$ by a function $\alpha'$ which acts only on $b$. Our definition is obviously inspired on the standard notion due to Kuratowski. In standard set theory $(a,b)$ is a set $\{\{a\}, \{a,b\}\}$ such that neither $a$ nor $b$ belong to $(a,b)$. In Flow, on the other hand, an ordered pair $(a,b)$ of the first kind is a function which does not act neither on $a$ nor on $b$.

Nevertheless, the second kind of ordered pair shows our approach is not equivalent to Kuratowki's. In the case where $\alpha = b$, we have the diagram to the right of {\sc Figure 2}. This non-Kuratowskian ordered pair $g = (a,b)$ acts on $b$, although it does not act on $a$. And no Kuratowskian ordered pair $(a,b)$ ever acts on either $a$ or $b$. So, in the general case, no ordered pair $(a,b)$ ever acts on $a$. That means:

\begin{teorema}\label{aa}

Any ordered pair $(a,a)$ is Kuratowskian.

\end{teorema}

The proof is straightforward. Since any ordered pair $(a,b)$ is a function, for the sake of abbreviation we write $x(a,b)$ for $x((a,b))$, for a given function $x$.

The reader can observe, from {\sc Figure 2}, that in a non-Kuratowskian ordered pair $f$, $\beta$ is the non-$\underline{\mathfrak{0}}$ $\mathfrak{F}$-successor of $\alpha$. In other words:

\begin{teorema}

$f = (a,b)$ is a non-Kuratowskian ordered pair iff $f$ acts only on a function $\alpha$ - which, in its turn, acts on one single term $a$ - and on its $\mathfrak{F}$-successor $\sigma_{\alpha}$, where $\sigma_{\alpha}\neq\underline{\mathfrak{0}}$.

\end{teorema}

\begin{description}
\item[Proof] If $f = (a,b)$, then $f$ acts at most on $\alpha$ and $\beta$, where $\alpha$ acts only on $a$ and $\beta$ acts only on $a$ and $b$. Suppose $f = (a,b)$ is  non-Kuratowskian. Then, $a\neq b$ (Theorem \ref{aa}); and $f = (a,\alpha)$, since $\alpha = b$. But $a\neq b$ entails $a\neq\alpha$. And since $\alpha$ acts on $a$, then $\beta$ acts on two terms: $\alpha$ and $a$. But that is the condition given for granting $\beta = \sigma_{\alpha}$, where $\sigma_{\alpha}\neq\underline{\mathfrak{0}}$. Finally, if $f$ acts only on $\alpha$ and $\sigma_{\alpha}$, where $\alpha$ acts only on $a$ and $\sigma_{\alpha}\neq\underline{\mathfrak{0}}$, then $f = (a,\alpha)$ is non-Kuratowskian.
\end{description}

Suppose $f = \underline{\mathfrak{1}}\big|_{x = \varphi_3 \vee x = \varphi_4}$. In that case $f$ acts on $\varphi_3$ and on its $\mathfrak{F}$-successor $\varphi_4$. Nevertheless, $\varphi_3$ does not act on just one single term. Thus, such an $f$ is not an ordered pair, let alone a non-Kuratowskian ordered pair.

\begin{teorema}

$(a,b) = (c,d)$ iff $a = c$ and $b = d$.

\end{teorema}

\begin{description}
\item[Proof] Straightforward from Definition \ref{orderedpair} and {\sc F9}$_{\mathbb{E}}$.
\end{description}

\begin{teorema}\label{teoreminhaparordenado}
If $a$ and $b$ are both emergent, then there is $f = (a,b)$.
\end{teorema}

\begin{description}
\item[Proof] From {\sc F9}$_{\mathbb{E}}$ we define the proper restriction $\beta$ of $\underline{\mathfrak{1}}$ for ``$x = a \vee x = b$'' as formula $F(x)$. So, $\beta (a) = a$, $\beta (b) = b$, $\beta (\beta) = \beta$, and $\beta (x) = \underline{\mathfrak{0}}$ for all remaining values of $x$. Analogously, the proper restriction $\alpha$ of $\underline{\mathfrak{1}}$ for ``$x = a$''gives us $\alpha (a) = a$, $\alpha (\alpha) = \alpha$, and $\alpha (x) = \underline{\mathfrak{0}}$ for the remaining values of $x$. And, from {\sc F9}$_{\mathbb{E}}$, we have $\alpha\neq a$. Analogously, we have $\beta\neq a$ and $\beta\neq b$. Finally, the proper restriction $f$ of $\underline{\mathfrak{1}}$ for ``$x = \alpha \vee x = \beta$ gives us $f(\alpha) = \alpha$, $f(\beta) = \beta$, $f(f) = f$, and $f(x) = \underline{\mathfrak{0}}$ for all the remaining values of $x$. Besides, $f\neq \alpha$ and $f\neq \beta$. But function $f$ is exactly that one in Definition \ref{orderedpair}. Hence, $f = (a,b)$.
\end{description}

To get an ordered pair $(\underline{\mathfrak{0}},\underline{\mathfrak{0}})$, all we have to do is to consider $\varphi_1$, which acts only on $\varphi_0$. So, $\varphi_1=_{def}(\underline{\mathfrak{0}},\underline{\mathfrak{0}})$. If $f$ acts only on $\alpha$ and $\varphi_0$, where $\alpha$ acts only on $a\neq\underline{\mathfrak{0}}$, then we adopt the convention $f = (a,\underline{\mathfrak{0}})$. For example, $\varphi_2$ is the ordered pair $(\varphi_0,\underline{\mathfrak{0}})$. There are no ordered pairs of the form $(\underline{\mathfrak{0}},b)$, where $b\neq\underline{\mathfrak{0}}$.

\begin{teorema}\label{potenciaZF}

If $f$ is a ZF-set, then $\wp(f)$ is a ZF-set.

\end{teorema}

\begin{description}
\item[Proof] If $f$ is a ZF-set, then it is emergent. So, $f$ acts only on emergent functions and {\sc F8} says any $t$ who lurks $f$ is emergent. Now, let $h = \underline{\mathfrak{1}}\big|_{t\trianglelefteq f}$. Then, $h[t]\Leftrightarrow t\trianglelefteq f$. So, $h = \mathfrak{p}(f)$ (Definition \ref{definindomaximapotencia}). But {\sc F8} also says such an $h$ is emergent. From Definition \ref{definindowp}, $p = \wp(f) = \underline{\mathfrak{1}}\big|_{t\subseteq f}$. But $p$ lurks $h$. So, once again from {\sc F8}, $\sigma_p\neq\underline{\mathfrak{0}}$.
\end{description}

\begin{definicao}\label{definindouniaoarbitraria}

Let $\mathbb{Z}(f)$. The {\em arbitrary ZF-union of all terms $g$ where $f$ acts\/} (or {\em arbitrary union of $f$\/}, for short) is defined as $u = \bigcup_{f[g]} g = \underline{\mathfrak{1}}\big|_{\exists g(f[g]\wedge g[t])}$.

\end{definicao}

The idea of arbitrary union is as it follows. If $f$ acts on any $g$ which acts on any $x$, then $u$ acts on that very same $x$; and if no $g$ acts on a given $x$, then $u$ does not act on that $x$. Particularly, we write $u = g\cup h$ for the case where $f$ acts at most on $g$ and $h$.

\begin{definicao}\label{intersecaoarbitraria}

Let $\mathbb{Z}(f)$. Then
$\bigcap_{f[g]}g =_{def} \underline{\mathfrak{1}} \big|_{\forall g(f[g] \Rightarrow g[t])}.$

\end{definicao}

Observe the definition above is equivalent to $\bigcap_{f[g]}g = \left( \bigcup_{f[g]} g\right) \big|_{\forall g(f[g] \Rightarrow g[t])}$. If $f$ acts only on $g$ and $h$, the arbitrary intersection may be written as $g\cap h$. In particular, for any $m$ and $n$ from language $\mathcal L$, $\varphi_m\cap\varphi_n = \varphi_m\circ\varphi_n$.

\begin{teorema}

$\bigcup_{\underline{\mathfrak{0}}[g]} = \bigcup_{\varphi_0[g]} = \bigcup_{\varphi_1[g]} = \varphi_0$.

\end{teorema}

Proof is straightforward. Last identity from last theorem illustrates our previous claim that it is possible the union $u$ be one of the terms where $f$ does act: $\varphi_1$ acts on $\varphi_0$, and $\bigcup_{\varphi_1[g]} = \varphi_0$.

\begin{teorema}\label{uniaoarbitrariadezfsets}

For any $f$, $\mathbb{Z}(f)\Rightarrow \mathbb{Z}(\bigcup_{f[g]} g)$.

\end{teorema}

\begin{description}
\item[Proof] If $u = \bigcup_{f[g]} g = \underline{\mathfrak{1}}\big|_{\exists g(f[g]\wedge g[t])}$, then it satisfies last conjunction of {\sc F8}. Thus, $\sigma_u\neq\underline{\mathfrak{0}}$. That entails $u$ is a ZF-set.
\end{description}

\begin{teorema}\label{intersecaoordinaria}

Let $x$ and $y$ be ZF-sets. Then, $x\cap y$ acts on $t$ iff $x[t]\wedge y[t]$.

\end{teorema}

The proof is immediate. Besides, if $x$ and $y$ are ZF-sets, then $\mathbb{Z}(x\cap y)$.

\begin{teorema}

For any pair $m$ and $n$ from language $\mathcal L$, $\varphi_m\cup\varphi_n = \varphi_{\mbox{\rm max}\{m,n\}}$, where $\mbox{\rm max}\{m,n\}$ denotes the maximum value between $m$ and $n$ relatively to $\preceq$.

\end{teorema}

\begin{description}
\item[Proof] Straightforward from Definition \ref{definindouniaoarbitraria}.
\end{description}

\begin{teorema}

For any $g$ and $h$ we have $g\cup h = h\cup g$.

\end{teorema}

\begin{description}
\item[Proof] Straightforward from Definition \ref{definindouniaoarbitraria}.
\end{description}

Observe for any $f\neq\underline{\mathfrak{0}}$ we have $f\cup \varphi_0 = f$.

Let $\alpha(x,y)$ be a formula where there is at least one occurrence of $x$, one occurrence of $y$, and all of them are free. Then, next formula is an axiom.

\begin{description}

\item[\sc F10$_{\alpha}$ - Creation] $\forall x\exists! y(\alpha(x,y))\Rightarrow \forall f(\mathbb{E}(f)\Rightarrow \exists g\forall x\forall y((\alpha(x,y)\wedge x\neq g)\Rightarrow ((f[x]\Rightarrow g(x) = y)\wedge (\neg f[x]\Rightarrow g(x) = \underline{\mathfrak{0}})))\wedge \sigma_g \neq \underline{\mathfrak{0}})$.

\end{description}

This last postulate is supposed to define functions $g$ from certain formulas $\alpha$ and any $f$, as long $f$ is emergent. We say $g$ {\em is created by $f$ and $\alpha(x,y)$\/}, and denote this by $g = f\big|^{\alpha(x,y)}$ or, simply, $g = f\big|^{\alpha}$.
Postulate {\sc F10}$_{\alpha}$ motivates us to employ an alternative notation for functions. When we see fit to do so, we may eventually write $a\marcio{f}b$ meaning $f(a) = b$. One obvious advantage is that it can be used to easily write down strings like $a_1\marcio{f} a_2\marcio{f} a_3\marcio{f} \cdot\cdot\cdot$, meaning $f(a_1) = a_2$, $f(a_2) = a_3$, and so on.

{\sc F10}$_{\alpha}$ allows us to get restrictions by other means besides {\sc F9}$_{\mathbb{E}}$. For example, from {\sc F10}$_{\alpha}$ there is a $g$ such that $\varphi_1\marcio{g} \varphi_2\marcio{g} \varphi_3\marcio{g} \underline{\mathfrak{0}}$. In that case $g\subset\sigma$. Function $g$ could be created, e.g., from $\varphi_4$.

With the aid of {\sc F10}$_{\alpha}$ it is quite easy to show $\mathfrak{F}$-composition is not commutative. As an example, consider $f$ as a function given by $\varphi_1\marcio{f}\varphi_2\marcio{f}\varphi_2$ and $g$ given by $\varphi_2\marcio{g}\varphi_3\marcio{g}\varphi_3$. In that case $g\circ f$ is given by $\varphi_1\marcio{g\circ f}\varphi_3\marcio{g\circ f}\underline{\mathfrak{0}}$ and $\varphi_2\marcio{g\circ f}\varphi_3\marcio{g\circ f}\underline{\mathfrak{0}}$, while $f\circ g$ is simply $\varphi_0$ (which is obviously different of $g\circ f$).

However, there are some functions which cannot exist, even with the comprehensive character of {\sc F10}$_{\alpha}$, as we can see in the next theorem.

\begin{teorema}\label{proibindoumbarracomt}

There is no $x$ where, for a given $t\neq\underline{\mathfrak{1}}$, $t\marcio{x}\underline{\mathfrak{1}}\marcio{x} t$.

\end{teorema}

\begin{description}
\item[Proof] Suppose there is such a function $x$. Then, from {\sc F5}, $x\circ x$ is given by $t\marcio{x\circ x} t$ and $\underline{\mathfrak{1}}\marcio{x\circ x} \underline{\mathfrak{1}}$. But from Theorem \ref{teoremamarcio2} we have that any function $h$ such that $h(\underline{\mathfrak{1}}) = \underline{\mathfrak{1}}$ entails $h = \underline{\mathfrak{1}}$. Thus, $x\circ x$ is supposed to be $\underline{\mathfrak{1}}$. But since $x$ is different of $\underline{\mathfrak{1}}$, then $x\circ x$ cannot be $\underline{\mathfrak{1}}$, according to Theorem \ref{dezenovedemarco}. That is a contradiction!
\end{description}

\begin{definicao}\label{localmenteinjetiva}

$f$ is {\em injective\/} iff $\forall r \forall s ((f[r]\wedge f[s] \wedge r\neq s )\Rightarrow f(r)\neq f(s))$. We denote this by $\mathbb{I}(f)$.

\end{definicao}

For example, any $\varphi_n$ is injective.

\begin{definicao}

If $\mathbb{E}(f)$, $Dom_f^{\mathfrak{F}} = \underline{\mathfrak{1}}\big|_{f[t]}$ and $Im_f^{\mathfrak{F}} = \underline{\mathfrak{1}}\big|_{\exists x(f(x) = t)}$ are, respectively, the {\em $\mathfrak{F}$-domain\/} and the {\em $\mathfrak{F}$-image\/} of $f$. We denote $f$ as $f:Dom_f^{\mathfrak{F}}\to Im_f^{\mathfrak{F}}$.

\end{definicao}

When we write $f:x\to y$, that means $x = Dom_f^{\mathfrak{F}}$ and $y = Im_f^{\mathfrak{F}}$, as long $f$ is emergent. If we recall that $\trianglelefteq$, $\mathbb{E}$, and $\mathbb{I}$ are given, respectively, in Definitions \ref{lurking}, \ref{definindoemergente}, and \ref{localmenteinjetiva}, next formula is our last axiom of this Section.

\begin{description}

\item[\sc F11 - $\mathfrak{F}$-Choice] $\forall f((\mathbb{E}(f)\wedge \exists r \exists s(r\neq s\wedge f[r]\wedge f[s]\wedge f(r)\neq f(s)))\Rightarrow \\ \exists c(c\trianglelefteq f\wedge \mathbb{I}(c)\wedge Dom_c^{\mathfrak{F}}\subseteq Dom_f^{\mathfrak{F}}\wedge Im_c^{\mathfrak{F}} = Im_f^{\mathfrak{F}}\wedge c\not\subseteq f))$.

\end{description}

This last postulate refers to emergent functions $f$ which act on more than one single term and such that there are at least two terms $r$ and $s$ where $f$ acts and $f(r)\neq f(s)$. In that case it says there always is a {\em choice-function\/} $c$ whose $\mathfrak{F}$-domain is a restriction of $f$'s $\mathfrak{F}$-domain but such that $c$ is injective. Besides, $c$ lurks $f$ (Definition \ref{lurking}) but it is no restriction of $f$. For example, let $f$ be given by $f(\varphi_0) = \varphi_7$, $f(\varphi_1) = \varphi_8$, and $f(\varphi_2) = \varphi_8$. In that case there may be an $c$ such that, e.g., $c(\varphi_0) = \varphi_8$ and $c(\varphi_1) = \varphi_7$.

Last axiom grants the existence of certain injective functions from arbitrary surjective functions whose $\mathfrak{F}$-images act on at least two terms.

\begin{description}

\item[\sc F11$_T$ - $\mathfrak{F}$-Trivial Choice] $\forall f((\mathbb{E}(f)\wedge \forall r \forall s((r\neq s\wedge f[r]\wedge f[s])\Rightarrow f(r)= f(s)))\Rightarrow \exists c(c\trianglelefteq f\wedge \mathbb{I}(c)\wedge Dom_c^{\mathfrak{F}}\subseteq Dom_f^{\mathfrak{F}}\wedge Im_c^{\mathfrak{F}} = Im_f^{\mathfrak{F}}))$.

\end{description}

This last formula is a complement for {\sc F11} in the sense it refers to the trivial case where all images $f(r)$ have the same constant value, if $f[r]$. In that case the choice $c$ is necessarily a restriction of $f$ such that $c$ acts on one single term.

\begin{teorema}[Partition Principle]\label{teoremadeparticao}

If $f:Dom_f^{\mathfrak{F}}\to Im_f^{\mathfrak{F}}$ is emergent, there is an injection $g:Im_f^{\mathfrak{F}}\to Im_g^{\mathfrak{F}}$, where $Im_g^{\mathfrak{F}}\subseteq Dom_f^{\mathfrak{F}}$.

\end{teorema}

\begin{description}
\item[Proof] First we consider the case $\exists r \exists s(f[r]\wedge f[s]\wedge f(r)\neq f(s))$. {\sc F11} says there is an injective choice-function $c:Dom_c^{\mathfrak{F}}\to Im_f^{\mathfrak{F}}$, where $Dom_c^{\mathfrak{F}}\subseteq Dom_f^{\mathfrak{F}}$. Since all terms involved are emergent, we can use {\sc F10}$_{\alpha}$ to define $g = u\big|^{\alpha(x,y)}$, where $u = Dom_c^{\mathfrak{F}}\cup Im_f^{\mathfrak{F}}$ and $\alpha(x,y)$ is the next formula: $(Im_c^{\mathfrak{F}}[x]\Rightarrow \exists t(c(t) = x\wedge y = t))\wedge (\neg Im_c^{\mathfrak{F}}[x]\Rightarrow y = \underline{\mathfrak{0}})$. Formula $\alpha(x,y)$ satisfies the demands from {\sc F10}$_{\alpha}$, since $c$ is injective. And that fact entails $g$ is injective as well. Observe both $c$ and $g$ lurk $f$, although $c$ is no restriction of $f$. The case where $\forall r \forall s((r\neq s\wedge f[r]\wedge f[s])\Rightarrow f(r)= f(s))$ is analogous to the previous one. The only difference is that now $Dom_c^{\mathfrak{F}}$ and $Im_f^{\mathfrak{F}}$ act on one single term, each. Besides, in that case $c$ is necessarily a restriction of $f$.
\end{description}

Theorem \ref{teoremadeparticao} is the Flow-theoretic version of the well-known Partition Principle. Its correspondence to the Partition Principle in the sense of ZF depends on considerations to be discussed in Sections \ref{standard} and \ref{pipiopatica}.

\begin{definicao}\label{pertencercomoagir}
$x\in f$ iff $f\subset\underline{\mathfrak{1}}\wedge \mathbb{Z}(x)\wedge f[x]$.
\end{definicao}

The negation of formula $x\in f$ is abbreviated as $x\not\in f$. We read $x\in f$ as ``$x$ belongs to $f$'' or ``$x$ is a member of $f$''. The symbol $\in$ is called {\em membership relation\/}. In Section \ref{standard} we prove ZF is `immersed' within Flow. Thanks to last definition it is possible, at least in principle, to prove ZFU (ZF with {\em Urelemente\/}, or {\em atoms\/}) is immersed within Flow as well. After all, there can be many terms $x$ such that no one belongs to them as long those $x$ are not ZF-sets. A simple example is $\sigma$ (and all their non-trivial restrictions), which is no restriction of $\underline{\mathfrak{1}}$. But that is a task for future papers.

\begin{teorema}
{\sc (i)} $\forall x (x\not\in x)$; {\sc (ii)} $\forall x(\underline{\mathfrak{0}}\not\in x)$; {\sc (iii)} $\forall x (\underline{\mathfrak{1}}\not\in x)$.
\end{teorema}

\begin{description}
\item[Proof] Item {\sc (i)} is consequence from {\sc F2}, since no $x$ acts on itself. Item {\sc (ii)} is consequence from Theorem \ref{fechamentorenato}: no $x$ acts on $\underline{\mathfrak{0}}$. Item {\sc (iii)} is immediate, since $\underline{\mathfrak{1}}$ is not a ZF-set.
\end{description}

Item {\sc (i)} looks like a proof that Flow is well-founded. Nevertheless, the issue of regularity is a little more subtle than that, as we can see later on.

$\{f,g\}=_{def}\underline{\mathfrak{1}}\big|_{x=f\vee x = g}$, where $\mathbb{E}(f) \wedge \mathbb{E}(g)$. So, {\em a pair\/} $\{f,g\}$ is a restriction of $\underline{\mathfrak{1}}$ which acts only on $f$ and $g$, as long they are both emergent.

\begin{teorema}\label{provandopares}

Let $\mathbb{E}(f) \wedge \mathbb{E}(g)$. If $u = \{ f,g\}$ is a pair, then $\mathbb{E}(u)$.

\end{teorema}

\begin{description}
\item[Proof] Let $h = \varphi_2\big|^{\alpha(x,y)}$ be defined from {\sc F10}$_{\alpha}$ by formula $\alpha(x,y)$ given by $(x = \varphi_0\Rightarrow y = f)\wedge (x = \varphi_1 \Rightarrow y = g)\wedge ((x\neq \varphi_0 \wedge x\neq \varphi_1)\Rightarrow y = \underline{\mathfrak{0}})$. According to {\sc F10}$_{\alpha}$, $\sigma_h \neq \underline{\mathfrak{0}}$. And since $h$ acts on emergent functions, then $h$ is emergent. Now, let $u = \underline{\mathfrak{1}}\big|_{t = g(\varphi_0)\vee t = g(\varphi_1)}$. That entails $u$ lurks $h$. And axiom {\sc F8} grants $\sigma_u \neq \underline{\mathfrak{0}}$. Hence, $u$ is emergent.
\end{description}

In particular, any pair of ZF-sets is a ZF-set.

On the left side of the image below we find a Venn diagram of a set $f$ (in the sense of ZF), with its elements $g$, $h$, and $i$. On the right side there is a diagram of a ZF-set $f$ (in the sense of Flow) which acts on $g$, $h$, and $i$. One of the main differences between both pictorial representations is the absence of arrows in the Venn diagram. Sets are `static' objects, which do nothing besides being elements of other sets. Venn diagrams do not emphasize the full role of the Extensionality Axiom from ZF: a set is defined either by its elements or by those terms who {\em do not belong\/} to it as well. Within Flow, however, it is explicitly emphasized that anything which is not inside the rectangle of a diagram has an image $\underline{\mathfrak{0}}$. Thus, in a sense, Flow diagrams were always somehow implicit within Venn diagrams. At the end, standard extensional set theories, like ZFC, are simply particular cases of a general theory of functions, as we can formally check at Subsection \ref{venn}.

\begin{picture}(60,45)

\put(33,20){\circle{200}}

\put(13,35){$f$}

\put(23,22){$g$}

\put(40,22){$h$}

\put(32,7){$i$}

\put(153,5){\framebox(100,30)}

\put(160,25){$f$}

\put(172,10){$g$}

\put(170,16){$\boldmath \curvearrowleft$}

\put(200,16){$\boldmath \curvearrowleft$}

\put(230,16){$\boldmath \curvearrowleft$}

\put(202,10){$h$}

\put(232,10){$i$}

\end{picture}

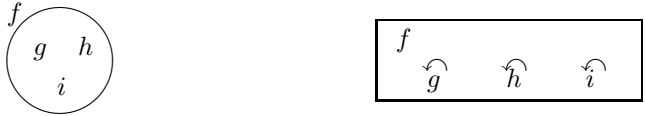
{\sc Figure 3:} Comparison between Venn diagrams and Flow diagrams.

\begin{teorema}

Any inductive function which acts only on ZF-sets is a ZF-set.

\end{teorema}

\begin{description}
\item[Proof] Immediate from the definitions of inductive function and ZF-set.
\end{description}

\begin{definicao}

$f$ is a {\em proper class\/} iff $\forall x(f[x]\Rightarrow \mathbb{Z}(x))\wedge \mathbb{C}(f)$.

\end{definicao}

In other words, no proper class is a ZF-set.

\begin{teorema}\label{unicidadevazio}
There is one single ZF-set $f$ such that for any $x$, we have $x\not\in f$.
\end{teorema}

\begin{description}
\item[Proof] $f = \varphi_0$. From {\sc F6}, $\varphi_0$ is unique. From Theorem \ref{noaction}, $\varphi_0$ is the only term different of $\underline{\mathfrak{0}}$ (recall $\underline{\mathfrak{0}}$ is not a ZF-set) which satisfies the condition given above. In other words, $\varphi_0$ is the {\em empty\/} ZF-set, which can be denoted by $\emptyset$.
\end{description}

\section{ZF is immersed in Flow}\label{standard}

We prove here vast portions of standard mathematics can be developed within Flow.

\subsection{ZFC axioms}

ZFC is a first-order theory with identity and with one predicate letter $f_1^2$, such that the formula $f_1^2(x,y)$ is abbreviated as $x\in y$, if $x$ and $y$ are terms, and is read as ``$x$ belongs to $y$'' or ``$x$ is an element of $y$''. The negation $\neg(x\in y)$ is abbreviated as $x\not\in y$. The axioms of ZFC are the following:

\begin{description}

\item[\sc ZF1 - Extensionality] $\forall x\forall y(\forall z(z\in x \Leftrightarrow z\in y)\Rightarrow x = y)$

\item[\sc ZF2 - Empty set] $\exists x\forall y(\neg(y\in x))$

\item[\sc ZF3 - Pair] $\forall x\forall y\exists z\forall t(t\in z\Leftrightarrow t = x \vee t = y)$

\end{description}

$x\subseteq y =_{def} \forall z(z\in x\Rightarrow z\in y)$

\begin{description}
\item[\sc ZF4 - Power set] $\forall x\exists y\forall z(z\in y\Leftrightarrow z\subseteq x)$
\end{description}

If $F(x)$ is a formula where there is no free occurrences of $y$, then:

\begin{description}
\item[\sc ZF5$_F$ - Separation] $\forall z\exists y\forall x(x\in y \Leftrightarrow x\in z\wedge F(x))$
\end{description}

The set $y$ is denoted by $\{x\in z/F(x)\}$.

If $\alpha(x,y)$ is a formula where all occurrences of $x$ and $y$ are free, then:

\begin{description}

\item[\sc ZF6$_\alpha$ - Replacement] $\forall x\exists!y\alpha(x,y)\Rightarrow\forall z\exists w\forall t(t\in w
\Leftrightarrow\exists s(s\in z\wedge\alpha(s,t)))$

\item[\sc ZF7 - Union set] $\forall x\exists y\forall z (z\in y\Leftrightarrow\exists t(z\in t\wedge t\in x))$

\end{description}

The set $y$ from {\sc ZF7} is abbreviated as $y = \bigcup_{t\in x} t$.

\begin{description}
\item[\sc ZF8 - Infinite] $\exists x(\emptyset\in x\wedge \forall y(y\in x\Rightarrow y\cup\{y\}\in x))$
\end{description}

\begin{description}
\item[\sc ZF9 - Choice] $\forall x(\forall y\forall z((y\in x\wedge z\in x\wedge y\neq z)\Rightarrow (y\neq\emptyset \wedge
y\cap z = \emptyset))\Rightarrow\\ \exists y\forall z(z\in x\Rightarrow \exists w (y\cap z = \{w\})))$
\end{description}

\subsection{ZF translation}\label{venn}

We refer to Flow as $\mbox{\boldmath{$\mathfrak{F}$}}$.

\begin{proposicao}\label{preservamatematica}
There is a translation from the language of ZF into the language of $\mbox{\boldmath{$\mathfrak{F}$}}$ such that every translated axiom of ZF is a theorem in $\mbox{\boldmath{$\mathfrak{F}$}}$.
\end{proposicao}

The translation from ZF (ZFC except Choice) into $\mbox{\boldmath{$\mathfrak{F}$}}$ is provided below.

\begin{center}
\begin{tabular}{|c|c|}\hline
\multicolumn{2}{|c|}{\sc Translating ZF into $\mbox{\boldmath{$\mathfrak{F}$}}$}\\ \hline \multicolumn{1}{|c|}{\sc ZF} & \multicolumn{1}{|c|}{$\mbox{\boldmath{$\mathfrak{F}$}}$}\\ \hline\hline
$\forall$ & $\forall_{\mathbb{Z}}$\\
$\exists$ & $\exists_{\mathbb{Z}}$\\
$x\in y$ & $y[x]$\\
$x\subseteq y$ & $x\subseteq y$\\ \hline
\end{tabular}
\end{center}

\noindent
where $\mathbb{Z}$ is the predicate ``to be a ZF-set'' from Definition \ref{defineZFconjunto}. Observe as well the translation of $\in$ follows Definition \ref{pertencercomoagir} due to the use of bounded quantifiers.

The proof of Proposition \ref{preservamatematica} is made through the following lemmas.

\begin{lema}\label{lemaprimeiro}

$\vdash_{\mbox{\boldmath{$\mathfrak{F}$}}} \mbox{``Translated}\; \mbox{{\sc ZF1}}''$.

\end{lema}

\begin{description}
\item[Proof] The translation of {\sc ZF1} is $\forall_{\mathbb{Z}} x\forall_{\mathbb{Z}} y(\forall_{\mathbb{Z}} z(x[z] \Leftrightarrow y[z])\Rightarrow x = y)$. If $x$ and $y$ are ZF-sets and $x[z]$ and $y[z]$, then $x(z) = z$ and $y(z) = z$ (Definition \ref{defineZFconjunto}). If $\neg x[z]$ or $\neg y[z]$, then either $x(z) = \underline{\mathfrak{0}}$ or $y(z) = \underline{\mathfrak{0}}$; or $z = x$ or $z = y$. So, the translated {\sc ZF1} considers the case where both $x$ and $y$ share the same images, except perhaps for $z = x$ or $z = y$. Thus, $x[z] \Leftrightarrow y[z]$ is equivalent to say that for any $z$ we have $x(z) = y(z)$, except perhaps for $z = x$ or $z = y$. That means $x\sim y$ (Definition \ref{definesim}). But Theorem \ref{marcionoslibertou} says $\underline{\mathfrak{0}}$ and $\varphi_0$ are the only clones. Since $\underline{\mathfrak{0}}$ is not a ZF-set, then the translation of {\sc ZF1} is a theorem for every ZF-set.
\end{description}

\begin{lema}

$\vdash_{\mbox{\boldmath{$\mathfrak{F}$}}} \mbox{``Translated}\; \mbox{\sc ZF2}''$.

\end{lema}

\begin{description}
\item[Proof] Translated {\sc ZF2} is $\exists_{\mathbb{Z}} x\forall_{\mathbb{Z}} y(\neg(x[y]))$: a corollary from Theorem \ref{unicidadevazio}.
\end{description}

\begin{lema}\label{lemadoparordenado}

$\vdash_{\mbox{\boldmath{$\mathfrak{F}$}}} \mbox{``Translated}\; \mbox{\sc ZF3}''$.

\end{lema}

\begin{description}
\item[Proof] Translated {\sc ZF3} is $\forall_{\mathbb{Z}} x\forall_{\mathbb{Z}} y\exists_{\mathbb{Z}} z\forall_{\mathbb{Z}} t(z[t]\Leftrightarrow (t = x \vee t = y))$. That is a consequence from Theorem \ref{provandopares}.
\end{description}

\begin{lema}

$\vdash_{\mbox{\boldmath{$\mathfrak{F}$}}} \mbox{``Translated}\; \mbox{\sc ZF4}''$.

\end{lema}

\begin{description}
\item[Proof] ]Translated {\sc ZF4} is $\forall_{\mathbb{Z}} x\exists_{\mathbb{Z}} y\forall_{\mathbb{Z}} z(y[z]\Leftrightarrow z\subseteq x)$. Theorem \ref{potenciaZF}.
\end{description}

\begin{lema}

$\vdash_{\mbox{\boldmath{$\mathfrak{F}$}}} \mbox{``Translated}\; \mbox{\sc ZF5}''$.

\end{lema}

\begin{description}
\item[Proof] Translated {\sc ZF5} is $\forall_{\mathbb{Z}} f\exists_{\mathbb{Z}} g\forall_{\mathbb{Z}} x(g[x] \Leftrightarrow f[x]\wedge F(x))$. If $g = f\big|_{F}$ ({\sc F9}$_{\mathbb{E}}$), then $g\subseteq f$. Besides, any restriction of a ZF-set is a ZF-set. Translated {\sc ZF5} demands no occurrence of $g$ in $F$. Recalling Observation \ref{consideracoessobrerestricao}, the calculation of $g = f\big|_{F}$ does not take into account any $t = g$. So, any occurrence of $g$ in $F$ becomes irrelevant.
\end{description}

\begin{lema}

$\vdash_{\mbox{\boldmath{$\mathfrak{F}$}}} \mbox{``Translated}\; \mbox{\sc ZF6}''$.

\end{lema}

\begin{description}
\item[Proof] Translated {\sc ZF6} is $\forall_{\mathbb{Z}} x\exists_{\mathbb{Z}}! y\alpha(x,y)\Rightarrow \forall_{\mathbb{Z}} z\exists_{\mathbb{Z}} w\forall_{\mathbb{Z}} t(w[t]
\Leftrightarrow\exists_{\mathbb{Z}} s(z[s]\wedge\alpha(s,t)))$ (we are supposed to rewrite {\sc ZF6} by means of $\forall$ and $\exists$, and only then translate it!). Let $g = z\big|^{\alpha(s,t)}$ be defined from {\sc F10}$_{\alpha}$ by formula $\alpha(s,t)$ provided by translated {\sc ZF6}. According to {\sc F10}$_{\alpha}$, $\sigma_g \neq \underline{\mathfrak{0}}$; and since $g$ acts on emergent functions, then $g$ is emergent. Now, let $w = \underline{\mathfrak{1}}\big|_{\exists s (z[s]\wedge t = g(s))}$. That entails $w$ lurks $g$. And axiom {\sc F8} grants $\sigma_w \neq \underline{\mathfrak{0}}$. Hence, $w$ is emergent. But $w$ is a restriction of $\underline{\mathfrak{1}}$ with $\sigma_w \neq \underline{\mathfrak{0}}$ and who acts only on ZF-sets. Therefore, $w$ is a ZF-set.
\end{description}

\begin{lema}\label{lemadauniaoZF}

$\vdash_{\mbox{\boldmath{$\mathfrak{F}$}}} \mbox{``Translated}\; \mbox{\sc ZF7}''$.

\end{lema}

\begin{description}
\item[Proof]The translated {\sc F7} is the formula $\forall_{\mathbb{Z}} f\exists_{\mathbb{Z}} u\forall_{\mathbb{Z}} t (u[t]\Leftrightarrow\exists_{\mathbb{Z}} r(r[t]\wedge f[r]))$. Once again we changed the names of the original variables in order to facilitate its reading. But the translated {\sc ZF7} is exactly Theorem \ref{uniaoarbitrariadezfsets}.
\end{description}

\begin{lema}

$\vdash_{\mbox{\boldmath{$\mathfrak{F}$}}} \mbox{``Translated}\; \mbox{\sc ZF8}''$.

\end{lema}

\begin{description}
\item[Proof] Translated {\sc ZF8} is $\exists_{\mathbb{Z}} x(x[\emptyset]\wedge \forall_{\mathbb{Z}} y(x[y]\Rightarrow x[y\cup\{y\}]))$. {\sc F7} states $\exists i ((\forall t (i(t) = t \vee i(t) = \underline{\mathfrak{0}}))\wedge \sigma_i \neq \underline{\mathfrak{0}} \wedge (i(\sigma_{\underline{\mathfrak{0}}}) = \sigma_{\underline{\mathfrak{0}}} \wedge \forall x (i(x) = x \Rightarrow (i(\sigma_x) = \sigma_x \neq\underline{\mathfrak{0}}))))$. Well, $\emptyset$ is exactly $\varphi_0$. So, $i[\varphi_0]$ and $i(\varphi_0) = \varphi_0$. Besides, $\mathbb{Z}(\varphi_0)$, and if $\mathbb{Z}(y)$ then $\sigma_y\neq\underline{\mathfrak{0}}$. And $y\cup \{y\}$ is exactly $\sigma_y$, where $\{y\} = \underline{\mathfrak{1}}\big|_{t = y}$. And the union of ZF-sets is a ZF-set (previous lemma). So, if $i$ acts on a ZF-set $t$, then $i$ acts on the ZF-set $\sigma_t$, which makes $i$ itself a ZF-set, since $\sigma_i\neq\underline{\mathfrak{0}}$.
\end{description}

\subsection{Grothendieck Universe and cardinals}\label{falandonisso}

Next we follow a maneuver somehow inspired on Cantor-Schr\"oder-Bernstein theorem \cite{Kolmogorov-75}.

\begin{definicao}\label{equipotencia}

$g$ and $h$ are {\em equipotent\/}, denoted by $g\equiv h$, iff there is $\tau$ where, for any $t$: {\sc (i)} $\tau[t]\Rightarrow\tau(\tau(t)) = t$; {\sc (ii)} $g[t]\Rightarrow (\tau[t] \wedge h[\tau(t)])$; {\sc (iii)} $h[t]\Rightarrow (\tau[t] \wedge g[\tau(t)])$; {\sc (iv)} $\tau[t]\Rightarrow (g[t]\vee h[t])$. We call function $\tau$ a {\em connector of $g/h$\/}.

\end{definicao}

For example, any $\varphi_n$ is equipotent to itself. And a connector of $\varphi_n/\varphi_n$ is $\varphi_n$.

\begin{definicao}\label{funcaoinfinita}

$g$ is {\em finite\/} iff $\exists h(g\equiv h \wedge \tau\not\equiv g)$, where $\tau$ is a connector of $g/h$ (Definition \ref{equipotencia}); and $g$ is {\em infinite\/} iff it is not finite.

\end{definicao}

As an example, $h = \varphi_5\big|_{t = \varphi_3\vee t = \varphi_4}$ acts only on $\varphi_3$ and $\varphi_4$. So, if we prove the existence of $\tau$ such that $\tau(\varphi_0) = \varphi_3$, $\tau(\varphi_3) = \varphi_0$, $\tau(\varphi_1) = \varphi_4$, and $\tau(\varphi_4) = \varphi_1$, then $h\equiv \varphi_2$ with connector $\tau$. Besides, $\tau\not\equiv\varphi_2 \wedge \tau\not\equiv h$, since $\tau$ acts on four terms, while $\varphi_2$ and $h$ act on two terms each. To prove the existence of such a $\tau$ we use {\sc F10}$_{\alpha}$, some pages below. But we hope the intuition is insightful.

\begin{teorema}

$\equiv$ is reflexive among emergent functions.

\end{teorema}

\begin{description}
\item[Proof] If $f$ is $\varphi_0$, assume $\underline{\mathfrak{0}}$ as connector, according to Definition \ref{equipotencia}. For the remaining cases, make $\tau = \underline{\mathfrak{1}}\big|_{f[x]}$ as a connector of $f/f$.
\end{description}

\begin{teorema}\label{equipotenciaquelevaaequipotencia}

$\forall f\forall g ((\mathbb{Z}(f)\wedge \mathbb{Z}(g)\wedge f\equiv g)\Rightarrow \sigma_f\equiv \sigma_g)$.

\end{teorema}

\begin{description}
\item[Proof] If $f$ and $g$ are ZF-sets, then $f\cup g$ is a ZF-set (Theorem \ref{uniaoarbitrariadezfsets}). That means $\sigma_f\neq\underline{\mathfrak{0}}$, $\sigma_g\neq\underline{\mathfrak{0}}$, and $\sigma_{f\cup g}\neq\underline{\mathfrak{0}}$. If $f\equiv g$, then there is a connector $\tau$ of $f/g$. Since $f\cup g$ is uncomprehensive, we can apply axiom {\sc F10}$_{\alpha}$ over $f\cup g$, where formula $\alpha(x,y)$ is $\tau(x) = y$. That entails, from {\sc F10}$_{\alpha}$, that $\sigma_{\tau}\neq\underline{\mathfrak{0}}$. Besides, $\sigma_f$, $\sigma_g$, and $\sigma_{f\cup g}$ are uncomprehensive and ZF-sets. That entails we can define a new connector $\tau'$ for $\sigma_f/\sigma_g$. All we have to do is to apply {\sc F10}$_{\alpha}$ over $\sigma_f\cup\sigma_g$ with formula $\alpha'(x,y)$ given by $\alpha(x,y)$ when $x\neq f$ and $x\neq g$, and such that $\alpha'(f,g)$ and $\alpha'(g,f)$. Hence, $\sigma_f\equiv\sigma_g$, where $\tau'$ is the connector of $\sigma_f/\sigma_g$.
\end{description}

The concept of equipotent functions is useful for different purposes. In this first paper we use it to discuss {\em cardinality\/}, as it follows in the next paragraphs. We mimic the construction originally introduced by von Neumann.

\begin{teorema}

$r = \underline{\mathfrak{1}}\big|_{\mathbb{Z}(t)\wedge \forall x(t[x]\Rightarrow x\subseteq t)\wedge (\exists x (t[x]\Rightarrow \exists y(t[y]\wedge\nexists s(t[s]\wedge y[s]))))}\Rightarrow\sigma_r = \underline{\mathfrak{0}}$.

\end{teorema}

\begin{description}
\item[Proof] Restriction $r$ above is a conjunction of three formulas. The first one, $\mathbb{Z}(t)$, grants the applicability of {\sc F9}$_{\mathbb{E}}$. Now, suppose $\sigma_r\neq\underline{\mathfrak{0}}$. Then $\mathbb{Z}(r)$, since $r$ acts only on ZF-sets. Observe $r$ acts on $\varphi_0$ by vacuity. And, for any $t$, $r[t]\Rightarrow r[\sigma_t]$, since $t\subseteq \sigma_t$ (where $\sigma_t\neq\underline{\mathfrak{0}}$). Thus, for any $x$, $r[x]\Rightarrow x\subseteq r$ ({\em transitivity\/}). That means $r$ satisfies the second conjunction of the formula used in the restriction above as well. Since $r$ acts on $\varphi_0$, then $r$ satisfies the third and final conjunction ({\em well-foundedness\/}). After all, $\varphi_0$ does not act on any term. So, if $F(t)$ is the formula used in the restriction above, then $F(r)$. But according to item ({\sc ii}) of Definition \ref{delimitadorarestricao}, that entails $\underline{\mathfrak{1}}(\sigma_r) = \underline{\mathfrak{0}}$. And such a condition is satisfied only for $\sigma_r = \underline{\mathfrak{0}}$, a contradiction. Therefore, $r$ cannot be a ZF-set.
\end{description}

\begin{definicao}\label{definicaodeordinais}

$\varpi = \underline{\mathfrak{1}}\big|_{\mathbb{Z}(t)\wedge \forall x(t[x]\Rightarrow x\subseteq t)\wedge (\exists x (t[x]\Rightarrow \exists y(t[y]\wedge\nexists s(t[s]\wedge y[s]))))}$ is the {\em ordinal function\/}. If $\varpi(t)\neq\underline{\mathfrak{0}}$, we say $t$ is an {\em ordinal\/}.

\end{definicao}

\begin{teorema}\label{ordinaisbemordenados}

Ordinal function $\varpi$ is well-ordered with respect to actions.

\end{teorema}

\begin{description}
\item[Proof] Suppose $\varpi[u]$ and $\varpi[v]$, where $u\neq v$. From the definition of $\varpi$ it follows that either $u\subseteq v$ or $v\subseteq u$. If $u[v]$ and $v[u]$, then $u(v) = v$ and $v(u) = u$, since $\varpi$ acts only on restrictions of $\underline{\mathfrak{1}}$. But that would entail $u = v$, according to {\sc F1}. Thus, if $\varpi[u]\wedge \varpi[v]$, then either $u[v]$, or $v[u]$ or $u = v$; but it is never the case $u[v]$ and $v[u]$. That entails actions define a total order if we are constraining ourselves to all terms where $\varpi$ acts. Besides, such a total order satisfies trichotomy. Concerning transitivity, if $\varpi$ acts on $u$, $v$, and $w$ in a way such that $u\subseteq v$ and $v\subseteq w$, the definition of restriction grants $u\subseteq w$; and the definition of $\varpi$ grants that $(v[u]\wedge w[v])\Rightarrow w[u]$. Finally, well-foundedness is granted from the definition of $\varpi$, since $\varpi[t]\Rightarrow \exists x (t[x]\wedge \nexists y (t[y]\wedge x[y]))$.
\end{description}

If the ordinal function $\varpi$ acts on any $t$, then that $t$ admits a least element with respect to actions. Besides, $\varpi$ satisfies the least-upper-bound property, i.e., any restriction of $\varpi$ with an upper bound has a least upper bound.

\begin{definicao}\label{definindocardinalidade}

Function $\mathfrak{c}$ is the {\em cardinality\/} of $g$, and we denote this by $\mathfrak{c} = |g|$, iff $\mathfrak{c}$ is the least term (with respect to action) such that $\varpi[\mathfrak{c}]$ and $\mathfrak{c}\equiv g$. A function $\mathfrak{c}$ is a {\em cardinal\/} iff it is the cardinality of some function $g$.

\end{definicao}

$|\varphi_3| = \varphi_3$. As another example, consider $\omega = \bigcap_{g\subseteq i\wedge \mbox{g is inductive}} g$ (Definition \ref{intersecaoarbitraria}), where $i$ is an arbitrary inductive (Definition \ref{funcaoindutiva}) ZF-set. If we use analogous procedures from standard set theories, we can easily identify $\omega$ with the ZF-set of natural numbers. And once again we have $|\omega| = \omega$, since $\varpi$ acts on $\omega$ and there is no $g$ such that $\sigma_g = \omega$. That means $\omega$ is an inductive function and, therefore, $\sigma_{\omega}\neq\underline{\mathfrak{0}}$. Hence, $\omega$ is a ZF-set. The details are analogous to those usually employed in ZFC, if we recall that $\omega[n]$ can be rewritten as $n\in \omega$, and $\sigma_n$ as $n\cup \{ n\}$. We can name $|\omega|$ simply as $\aleph_0$.

\begin{teorema}

The cardinality of any ZF-set, if it exists, is a ZF-set.

\end{teorema}

\begin{description}
\item[Proof] A sketch for the proof: If $f$ is {\em finite\/} (see Definition \ref{funcaoinfinita}), then $|f| = \varphi_n$ for some $n$. So, the cardinality of any finite $f$ is a ZF-set. If $f$ is equipotent to $\omega$, then $f$ is infinite and $|f| = \omega$, where $\omega$ is a ZF-set. According to Theorem \ref{potenciaZF}, the restricted power of any ZF-set is a ZF-set. Besides, the restricted power of any ordinal is an ordinal. If the ordinal is a ZF-set, then its restricted power is a ZF-set as well. Besides, if $f$ is a ZF-set and $|f| = \mathfrak{c}$, then the restricted power of $f$ is equipotent to the restricted power of $\mathfrak{c}$. Analogous rationale can be used for arbitrary unions, since the union of ZF-sets is a ZF-set (Theorem \ref{uniaoarbitrariadezfsets}), and the arbitrary union of ordinals is an ordinal (this last statement can be proved in an analogous way how ordinals are coped in standard literature \cite{Jech-03}).
\end{description}

A complete proof of last theorem demands a lot of room unavailable in a paper as long as this. We hope the sketch above is enough for persuading the reader.

\begin{teorema}

$\underline{\mathfrak{1}}\big|_{\mathbb{Z}(x)}$ exists and it has $\mathfrak{F}$-successor $\underline{\mathfrak{0}}$.

\end{teorema}

\begin{description}
\item[Proof] Let $u = \underline{\mathfrak{1}}\big|_{\mathbb{Z}(x)}$. We can use {\sc F9}$_{\mathbb{E}}$, since any ZF-set is emergent. Suppose $\sigma_u \neq \underline{\mathfrak{0}}$. Then $u$ is a ZF-set. Therefore, $\sigma_u$ acts on all ZF-sets. Thus, $\sigma_u = u$, which is a contradiction.
\end{description}

\begin{teorema}\label{universogrothen}

$\underline{\mathfrak{1}}\big|_{\mathbb{Z}(x)}$ is a Grothendieck Universe.

\end{teorema}

\begin{description}
\item[Proof] Let us denote $\underline{\mathfrak{1}}\big|_{\mathbb{Z}(x)}$ by $u$. Remember $r\in s$ iff $s\subset\underline{\mathfrak{1}}\wedge \mathbb{Z}(r)\wedge s[r]$. For any $x$, if $u[x]$, then $x$ is a ZF-set. And according to Definition \ref{defineZFconjunto}, if $x[y]$ then $y$ is a ZF-set as well. From Lemma \ref{lemadoparordenado} we know that if $u[x]$ and $u[y]$, then $u[\{x,y\}]$. From Theorem \ref{potenciaZF} we know that if $u[x]$, then $u[\wp(x)]$. Finally, according to Lemma \ref{lemadauniaoZF}, if $u$ acts on all elements of an arbitrary ZF-set $f$ (which is equivalent to the notion of a ZF-family of ZF-sets), then $u$ acts on $\bigcup_{f[g]} g$.
\end{description}

If the reader is puzzled by the meaning of $\underline{\mathfrak{1}}\big|_{\mathbb{Z}(x)}$, we are talking about Definition \ref{defineZFconjunto}. That means $\underline{\mathfrak{1}}\big|_{\mathbb{Z}(x)} = \underline{\mathfrak{1}}\big|_{\sigma_x\neq\underline{\mathfrak{0}}\wedge \forall t(x[t]\Rightarrow (x(t) = t \wedge \mathbb{Z}(t)))}$, where $\mathbb{Z}(t)$ is recursively defined by $\sigma_t\neq\underline{\mathfrak{0}}\wedge \forall r(t[r]\Rightarrow (t(r) = r \wedge \mathbb{Z}(r)))$. Nothing prevents us from assuming $\mathbb{Z}$ as a monadic predicate letter to be added to Flow's language. That would take place if we assumed the definition of $\mathbb{Z}$ as an ampliative one. But we prefer to assume that as an unnecessary maneuver.

Let $u' = \underline{\mathfrak{1}}\big|_{\mathbb{Z}(t)\wedge \exists x (t[x]\wedge \nexists y(t[y]\wedge x[y]))}$. By using the same techniques illustrated above we can prove $\sigma_{u'} = \underline{\mathfrak{0}}$. The reader should observe that such an $u'$ is also a Grothendieck universe, with an extra condition: regularity. So, we refer to $u'$ as the {\em well-founded Grothendieck Universe\/}. Observe $u'\subseteq u$, where $u$ is the same function from Theorem \ref{universogrothen}.

Since the cardinality $|u|$ of a Grothendieck universe $u$ is sup$_{u[t]}|t|$, then Theorem \ref{universogrothen} grants the existence of strongly inaccessible cardinals, provided we are able to determine cardinalities.

\begin{definicao}

A {\em model of ZF\/} is $\mathfrak{m} = \langle u,[],\eqcirc\rangle$, where: {\sc (i)} $u\subset \underline{\mathfrak{1}}$; {\sc (ii)} $\forall x (\mathbb{Z}(x)\Leftrightarrow u[x])$; {\sc (iii)} $p = u\big|_{\exists a\exists b(u[a]\wedge u[b]\wedge t = (a,b))}$; {\sc (iv)} $[] = p\big|_{\forall a\forall b(t = (a,b)\Leftrightarrow a[b])}$; {\sc (v)} $\eqcirc = p\big|_{\forall a\forall b(t = (a,b)\Leftrightarrow a = b)}$.

\end{definicao}

It is easy to check $\sigma_{p} = \sigma_{[]} = \sigma_{\eqcirc} = \underline{\mathfrak{0}}$, if $\sigma_u = \underline{\mathfrak{0}}$. A model of ZF depends only on the term $u$. Function $[]$ maps membership $\in$ in ZF, in the same way we did in Subsection \ref{venn}. Function $\eqcirc$ maps identity in ZF. Both $[]$ and $\eqcirc$ depend on who is $u$. If $u$ is our Grothendieck universe (Theorem \ref{universogrothen}), we have at least one model of ZF, since we proved above ZF is immersed within Flow.

\section{The Partition Principle}\label{pipiopatica}

The Partition Principle (PP), within standard set theories, states: if there is a surjection $f:x\to y$, then there is an injection $g:y\to x$. The Axiom of Choice (AC) entails PP. Nevertheless, until now it has endured as the oldest open problem in set theory whether PP implies AC. For a review about the subject see, e.g., \cite{Banaschewski-90} \cite{Higasikawa-95} \cite{Howard-16} \cite{Pelc-78}.

Although this Section is based on $\mathfrak{F}$-Choice Axiom {\sc F11} and {\sc F11}$_T$, all our claims regarding the Partition Principle and the Axiom of Choice are expressible within translated ZF (Section \ref{standard}). Any function $f$ in the sense of ZF is a set of ordered pairs. That translates into Flow as an $f\subset\underline{\mathfrak{1}}$ who acts on ordered pairs $(a,b)$ (Definition \ref{orderedpair}), where both $a$ and $b$ are ZF-sets. Besides, $f$ itself is a ZF-set. But {\sc F10}$_{\alpha}$ grants the existence of an emergent function $g$ such that $g(a) = b$ whenever $f[(a,b)]$. So, injective (Definition \ref{localmenteinjetiva}) emergent functions who act on ZF-sets, in the sense of Flow, correspond to bijective functions in the sense of ZF (recall we consider only $\mathfrak{F}$-images and not `codomains'). And those emergent functions who are not injective, in the sense of Flow, correspond to surjective functions who are not injective, in the sense of ZF.

Now, consider the next postulate:

\begin{description}

\item[Well-Foundedness] $\forall f((f\neq\underline{\mathfrak{0}}\wedge f\neq\varphi_0)\Rightarrow \exists x(f[x]\wedge \nexists t(f[t]\wedge x[t])))$.

\end{description}

This is analogous to the axiom of regularity (AR) in ZF, in the sense that the translation of AR into Flow's language is a straightforward theorem from our proposed Well-Foundedness Postulate and our Self-Reference Postulate: there can be no infinite chains of membership. If $f$ is a ZF-set, then there always is an $x$ such that $x$ belongs to $f$, but there is no intersection between $f$ and $x$.

\begin{definicao}

$\psi$ is a {\em hyperfunction\/}, and we denote this by $\mathbb{H}(\psi)$, iff: {\sc (i)} $\sigma_{\psi}\neq\underline{\mathfrak{0}}\wedge \psi\subset\underline{\mathfrak{1}} \wedge \exists t(\psi[t])$; {\sc (ii)} $\forall x(\psi[x]\Rightarrow (x\subset\underline{\mathfrak{1}}\wedge\psi[\sigma_x]\wedge \exists y(\psi[y]\wedge x[y])))$; {\sc (iii)} $\forall x(\psi[x]\Rightarrow \forall y ((y\subseteq x\wedge y\neq\varphi_0\wedge y\neq\underline{\mathfrak{0}})\Rightarrow \psi[y]))$.
\end{definicao}

Any $\psi$ which satisfies last definition is called a {\em non-well-founded\/} function as well. Item {(\sc ii}) entails a violation of Well-Foundedness.

Now, consider another possible postulate to be added to our formal system:

\begin{description}

\item[Hyperfunctions] $\exists \psi (\mathbb{H}(\psi)) \wedge \forall f(\mathbb{H}(f)\Rightarrow \mathbb{E}(f))\wedge \underline{\mathfrak{1}}\big|_{\forall h(\mathbb{H}(h)\Rightarrow h[t])} = \varphi_0$.

\end{description}

Obviously Well-Foundedness and Hyperfunctions are inconsistent with each other. Thus, we must choose which of them to consider, if we intend to add new postulates besides {\sc F1}$\sim${\sc F11}.

This last postulate says three things: {\sc (i)} there is at least one hyperfunction; {\sc (ii)} every hyperfunction is emergent; and {\sc (iii)} there is no common $t$ where all possible hyperfunctions act. We refer to this last item as the {\em Principle of Non-Commonality\/} among hyperfunctions. Clearly any hyperfunction $\psi$ acts on an infinity of terms: if $\psi$ acts on any $x$, then it acts on $\sigma_x$, as it happens in the Axiom of Infinity ({\sc F7}). Nevertheless, inductive functions (Definition \ref{funcaoindutiva} states inductive functions are those who satisfy {\sc F7}) do not follow any analogue of the Principle of Non-Commonality. If $\mathbb{J}$ is a predicate such that $\mathbb{J}(i)$ is equivalent to `$i$ is an inductive function', then $\underline{\mathfrak{1}}\big|_{\forall i(\mathbb{J}(i)\Rightarrow i[t])}$ is simply $\omega$, i.e., the unique function who acts on all $\varphi_n$ and only on them (the ZF-set of all finite ordinals). That happens because any inductive function $i$ acts on $\varphi_0$, $\varphi_1$, $\varphi_2$, and so on.

It is worth to remember that any emergent function $f\neq\varphi_0$ acts only on emergent functions. Thus, for the sake of illustration, let us show a few examples of hyperfunctions. For that purpose, it is easier to follow next Figure.


\begin{picture}(60,230)

\put(25,0){$\vdots$}

\put(25,90){$\vdots$}

\put(25,220){$\vdots$}

\put(75,0){$\vdots$}

\put(75,90){$\vdots$}

\put(75,220){$\vdots$}

\put(125,0){$\vdots$}

\put(125,90){$\vdots$}

\put(125,220){$\vdots$}

\put(185,0){$\vdots$}

\put(185,90){$\vdots$}

\put(185,220){$\vdots$}

\put(235,0){$\vdots$}

\put(235,90){$\vdots$}

\put(235,220){$\vdots$}

\put(285,0){$\vdots$}

\put(285,90){$\vdots$}

\put(285,220){$\vdots$}

\put(3,13){\framebox(45,35)}

\put(53,13){\framebox(45,35)}

\put(103,13){\framebox(45,35)}

\put(163,13){\framebox(45,35)}

\put(213,13){\framebox(45,35)}

\put(263,13){\framebox(45,35)}


\put(3,53){\framebox(45,35)}

\put(53,53){\framebox(45,35)}

\put(103,53){\framebox(45,35)}

\put(163,53){\framebox(45,35)}

\put(213,53){\framebox(45,35)}

\put(263,53){\framebox(45,35)}

\put(6,209){$a$}

\put(56,209){$b$}

\put(106,209){$c$}

\put(41,185){$b$}

\put(38,192){$\boldmath \curvearrowleft$}

\put(91,185){$c$}

\put(88,192){$\boldmath \curvearrowleft$}

\put(141,185){$a$}

\put(138,192){$\boldmath \curvearrowleft$}

\put(166,209){$j$}

\put(216,209){$k$}

\put(266,209){$l$}

\put(201,185){$k$}

\put(198,192){$\boldmath \curvearrowleft$}

\put(251,185){$l$}

\put(248,192){$\boldmath \curvearrowleft$}

\put(301,185){$j$}

\put(298,192){$\boldmath \curvearrowleft$}


\put(6,169){$\sigma_a$}

\put(6,105){$\sigma_a$}

\put(56,169){$\sigma_b$}

\put(56,105){$\sigma_b$}

\put(106,169){$\sigma_c$}

\put(106,105){$\sigma_c$}

\put(166,105){$\sigma_j$}

\put(216,105){$\sigma_k$}

\put(266,105){$\sigma_l$}

\put(41,145){$b$}

\put(38,152){$\boldmath \curvearrowleft$}

\put(91,145){$c$}

\put(88,152){$\boldmath \curvearrowleft$}

\put(141,145){$a$}

\put(138,152){$\boldmath \curvearrowleft$}

\put(25,145){$a$}

\put(23,152){$\boldmath \curvearrowleft$}

\put(75,145){$b$}

\put(73,152){$\boldmath \curvearrowleft$}

\put(125,145){$c$}

\put(123,152){$\boldmath \curvearrowleft$}


\put(166,169){$\sigma_j$}

\put(216,169){$\sigma_k$}

\put(266,169){$\sigma_l$}

\put(201,145){$k$}

\put(198,152){$\boldmath \curvearrowleft$}

\put(251,145){$l$}

\put(248,152){$\boldmath \curvearrowleft$}

\put(301,145){$j$}

\put(298,152){$\boldmath \curvearrowleft$}

\put(185,145){$j$}

\put(183,152){$\boldmath \curvearrowleft$}

\put(235,145){$k$}

\put(233,152){$\boldmath \curvearrowleft$}

\put(285,145){$l$}

\put(283,152){$\boldmath \curvearrowleft$}


\put(6,129){$\sigma_{\sigma_a}$}

\put(6,112){$\boldmath \curvearrowleft$}

\put(56,112){$\boldmath \curvearrowleft$}

\put(106,112){$\boldmath \curvearrowleft$}

\put(166,112){$\boldmath \curvearrowleft$}

\put(216,112){$\boldmath \curvearrowleft$}

\put(266,112){$\boldmath \curvearrowleft$}

\put(56,129){$\sigma_{\sigma_b}$}

\put(106,129){$\sigma_{\sigma_c}$}

\put(41,105){$b$}

\put(38,112){$\boldmath \curvearrowleft$}

\put(91,105){$c$}

\put(88,112){$\boldmath \curvearrowleft$}

\put(141,105){$a$}

\put(138,112){$\boldmath \curvearrowleft$}

\put(25,105){$a$}

\put(23,112){$\boldmath \curvearrowleft$}

\put(75,105){$b$}

\put(73,112){$\boldmath \curvearrowleft$}

\put(125,105){$c$}

\put(123,112){$\boldmath \curvearrowleft$}


\put(166,129){$\sigma_{\sigma_j}$}

\put(216,129){$\sigma_{\sigma_k}$}

\put(266,129){$\sigma_{\sigma_l}$}

\put(201,105){$k$}

\put(198,112){$\boldmath \curvearrowleft$}

\put(251,105){$l$}

\put(248,112){$\boldmath \curvearrowleft$}

\put(301,105){$j$}

\put(298,112){$\boldmath \curvearrowleft$}

\put(185,105){$j$}

\put(183,112){$\boldmath \curvearrowleft$}

\put(235,105){$k$}

\put(233,112){$\boldmath \curvearrowleft$}

\put(285,105){$l$}

\put(283,112){$\boldmath \curvearrowleft$}


\put(3,102){\framebox(45,35)}

\put(53,102){\framebox(45,35)}

\put(103,102){\framebox(45,35)}

\put(163,102){\framebox(45,35)}

\put(213,102){\framebox(45,35)}

\put(263,102){\framebox(45,35)}


\put(3,142){\framebox(45,35)}

\put(53,142){\framebox(45,35)}

\put(103,142){\framebox(45,35)}

\put(163,142){\framebox(45,35)}

\put(213,142){\framebox(45,35)}

\put(263,142){\framebox(45,35)}


\put(3,182){\framebox(45,35)}

\put(53,182){\framebox(45,35)}

\put(103,182){\framebox(45,35)}

\put(163,182){\framebox(45,35)}

\put(213,182){\framebox(45,35)}

\put(263,182){\framebox(45,35)}


\put(6,79){$d$}

\put(56,79){$e$}

\put(106,79){$f$}

\put(36,57){$\sigma_a$}

\put(36,62){$\boldmath \curvearrowleft$}

\put(36,17){$\sigma_a$}

\put(36,24){$\boldmath \curvearrowleft$}

\put(86,57){$\sigma_b$}

\put(86,62){$\boldmath \curvearrowleft$}

\put(86,17){$\sigma_b$}

\put(86,24){$\boldmath \curvearrowleft$}

\put(136,57){$\sigma_c$}

\put(136,62){$\boldmath \curvearrowleft$}

\put(136,17){$\sigma_c$}

\put(136,24){$\boldmath \curvearrowleft$}

\put(166,79){$m$}

\put(216,79){$n$}

\put(266,79){$o$}

\put(196,57){$\sigma_j$}

\put(196,62){$\boldmath \curvearrowleft$}

\put(196,17){$\sigma_j$}

\put(196,24){$\boldmath \curvearrowleft$}

\put(246,57){$\sigma_k$}

\put(246,62){$\boldmath \curvearrowleft$}

\put(246,17){$\sigma_k$}

\put(246,24){$\boldmath \curvearrowleft$}

\put(296,57){$\sigma_l$}

\put(296,62){$\boldmath \curvearrowleft$}

\put(296,17){$\sigma_l$}

\put(296,24){$\boldmath \curvearrowleft$}


\put(6,39){$g$}

\put(56,39){$h$}

\put(106,39){$i$}

\put(166,39){$p$}

\put(216,39){$q$}

\put(266,39){$r$}


\put(20,17){$a$}

\put(18,24){$\boldmath \curvearrowleft$}

\put(70,17){$b$}

\put(67,24){$\boldmath \curvearrowleft$}

\put(120,17){$c$}

\put(117,24){$\boldmath \curvearrowleft$}

\put(180,17){$j$}

\put(177,24){$\boldmath \curvearrowleft$}

\put(230,17){$k$}

\put(227,24){$\boldmath \curvearrowleft$}

\put(280,17){$l$}

\put(277,24){$\boldmath \curvearrowleft$}

\end{picture}

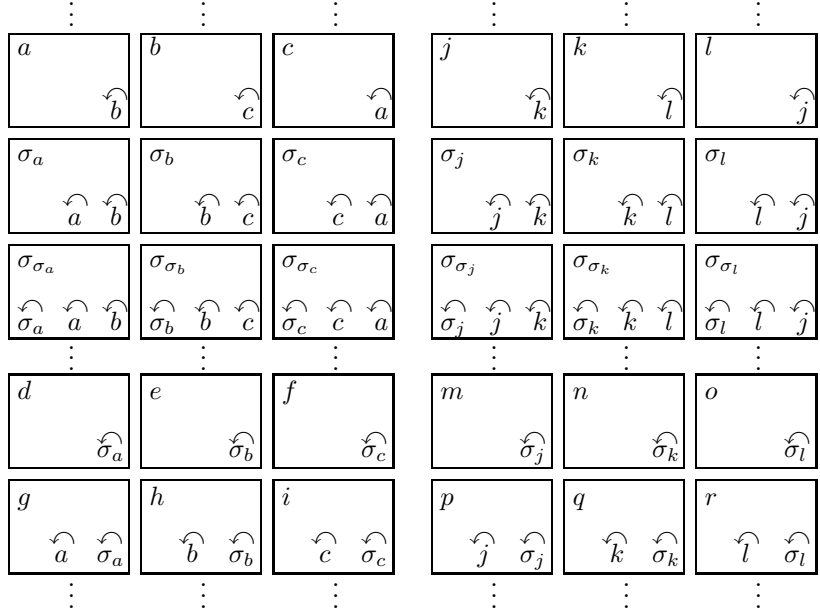
{\sc Figure 4:} Examples of hyperfunctions.\\

In the Figure above we suggest the possible existence of infinitely many terms where a given hyperfunction $\psi$ can act: $a$, $b$, $c$, $\sigma_a$, $\sigma_b$, $\sigma_c$, $\cdots$, $j$, $k$, $l$, $\sigma_j$, $\cdots$, $d$, $e$, $f$, $\cdots$. We represent this possibility by two `big columns', one to the right and one to the left, where each `big column' is formed by three `little columns'. Observe, for example, $d$ is a proper restriction of $\sigma_{\sigma_a}$. In such a case, $\psi$ is a hyperfunction, if for any $x$ we have $\psi[x]\Rightarrow \psi(x) = x$. Thus we have an infinite amount of infinite chains of membership relations: $\cdots b\in a\in c\in b\in a \in \cdots$; $\cdots b\in a \in c \in \sigma_c \in f \in \sigma_f \in \sigma_{\sigma_f}\in\cdots$; $\cdots$. That corresponds to the set-theoretic version of violation of the Axiom of Regularity. If we consider a function $\psi'$ who acts only on those terms related to $a$, $b$, and $c$ (the left `big column', where we find $a$, $b$, $c$, $\sigma_a$, $\sigma_b$, $\sigma_c$, $\sigma_{\sigma_a}$, $\sigma_{\sigma_b}$, $\sigma_{\sigma_c}$, $d$, $e$, $f$, etc.), $\psi'$ is a hyperfunction as well. The same happens to a function $\psi''$ who acts only on those terms related to $j$, $k$, and $l$ (the right `big column'). That means there can be proper restrictions of hyperfunctions who are hyperfunctions. On the other hand, each small column constitutes a proper restriction of $\psi$ as well. But none of such restrictions is a hyperfunction. In this simple example we started with ensembles (`big columns') of three terms ($a$, $b$, $c$; $j$, $k$, $l$). But we can obviously suggest other situations with as many terms as we wish, as long we keep in mind all properties a hyperfunction is supposed to hold. Nevertheless, there is no need to assume a hyperfunction is necessarily supposed to be like this. If we add a term $c'$ on the left `big column' who acts on both $a$ and $j$, we still have a hyperfunction $\psi$ which acts on all terms suggested above (besides $c'$). And, in that case, the left `big column' no longer represents a proper restriction of $\psi$ who is still a hyperfunction. Of course, such a distinction between left and right `big column' is just a pedagogical tool for our illustration purposes, since there is no clear criterium to distinguish them in the case $\psi$ acts on $c'$.

All of this helps us to understand the meaning of the Principle of Non-Commonality $\underline{\mathfrak{1}}\big|_{\forall h(\mathbb{H}(h)\Rightarrow h[t])} = \varphi_0$ introduced in axiom Hyperfunctions: there is no common term $x$ where all hyperfunctions act. That means our postulate Hyperfunctions has no constructive character, in the sense it does not provide any sufficiently clear criteria for building hyperfunctions. That postulate simply says there are terms like that. Our axiom of Weak Extensionality ({\sc F1}) still holds. But, in general, it is undecidable whether $\psi$ and $\upsilon$ are identical, in the case both are hyperfunctions. This fact grants the next theorem.

\begin{teorema}\label{teoremaqueantecedegrande}

Any hyperfunction $\psi$ admits an infinite proper restriction $\psi'\subset \psi$ which can be well-ordered.

\end{teorema}

\begin{description}
\item[Proof] According to Theorem \ref{teoremadeparticao} (see also the paragraph that follows the proof), if $f:\psi\to\varphi_1$ is a function (observe this is a `constant' function where all images $f(x)$ have the same value $\varphi_0$), then there is an injection $g:\varphi_1\to\psi$. Such an injection `chooses' a term $x_0$, according to {\sc F11}$_T$. Now we can recursively define a proper restriction $\psi'$ of $\psi$ as it follows: {\sc (i)} $\psi'(x_0) = x_0$; {\sc (ii)} $\forall x (\psi'[x]\Rightarrow \psi'(\sigma_x) = \sigma_x)$. Function $\sigma$ itself defines a well-order over $\psi'$.
\end{description}

Last theorem allows us to define a $f:\psi'\to\omega$ as it follows: {\sc (i)} $f(x_0) = \varphi_0$; {\sc (ii)} $\forall x\forall y(f(x) = y \Rightarrow f(\sigma_x) = \sigma_y)$. According to Theorem \ref{teoremadeparticao}, there is an injection $i:\omega\to\psi'$. But such an injection does not preserve the well-order from $\omega$, since the choice function $c:\psi'\to\omega$ obtained from {\sc F11} is no restriction of $f$. Besides, it is worth to remark such a proper restriction $\psi'$ of $\psi$ is not a hyperfunction.

Finally, we introduce the next postulate.

\begin{description}

\item[Hyper-ZF-Sets] $\forall f(\mathbb{H}(f)\Rightarrow \mathbb{Z}(f))$.

\end{description}

Any hyperfunction is a ZF-set (this axiom works together with Definition \ref{defineZFconjunto}). That means we have a new Grothendieck Universe $u = \underline{\mathfrak{1}}\big|_{\mathbb{Z}(t)}$ which works as a model of ZF without Regularity, without the Axiom of Choice of ZFC but with the Partition Principle. In such a model we prove next there are ZF-sets which cannot be well-ordered. Therefore, the translated Axiom of Choice from ZFC is not a theorem within our new model, termed $\mathfrak{F}'$ (based on axioms {\sc F1}$\sim${\sc F11}$+$Hyperfunctions$+$Hyper-ZF-Sets), while the Partition Principle is a consequence from our Axiom of $\mathfrak{F}$-Choice (Theorem \ref{teoremadeparticao}). We are not suggesting the Axiom of Choice in ZFC is somehow related to the Axiom of Regularity. Actually, it is not. But we chose to work within a non-well-founded version of ZF for two reasons: {\sc (i)} to introduce a non-well-founded model of ZF which cannot be modeled within a von Neumann Universe and, therefore, it is out of the scope of some metamathematical techniques like Scott's trick \cite{Scott-55} for determining cardinalities; and {\sc (ii)} hyperfunctions make it easier for us to prove our main result.

Next theorem makes use of the fact that AC in ZFC is equivalent to Zermelo's Theorem  (see, e.g., Theorem 5.1 on page 48 of \cite{Jech-03}).

\begin{teorema}\label{onzedeagostode2020}

There is at least one ZF-set in $\mathfrak{F}'$ that cannot be well-ordered.

\end{teorema}

\begin{description}
\item[Proof] Let $\psi$ be a hyperfunction. If $\psi$ can be well-ordered, then we can {\em define\/} an injection $f:\psi\to o$, where $o$ is an ordinal (see Definition \ref{definicaodeordinais} and Theorem \ref{ordinaisbemordenados}). Thus, the existence of such a function can be proven by transfinite induction in the sense of defining it by transfinite recursion. Suppose we have a surjection $f_0:\psi\to\varphi_1$. That can be defined by the next formula: for all $t$, $\psi[t]\Rightarrow f_0(t) = \varphi_0$. Theorem \ref{teoremadeparticao} (see discussion right after its proof as well) grants the existence of an injection $g_0$ which acts only on $\varphi_0$ (the only term where $\varphi_1$ acts), where $g_0(\varphi_0) = x_0$. Observe $x_0$ is one of the terms where $\psi$ acts. Observe as well that such an $g_0$ is the inverse of an injection $h_0$ which acts only on $x_0$ and such that $h_0(x_0) = \varphi_0$, according to {\sc F11}$_T$. It is rather important to observe that both {\sc F11}$_T$ and {\sc F11} do not simply `choose' injections, but also terms where $\psi$ acts. For the sake of argument (before we simply go to the case of all remaining ordinals), let us go to the next immediate step of our transfinite induction. Let $f_1:\psi\to \varphi_2$ be a surjective function defined as it follows: for all $t$, $(f_0[t]\Rightarrow f_1(t) = f_0(t)) \wedge ((\psi[t]\wedge \neg f_0[t])\Rightarrow f_1(t) = \varphi_1)$. In that case $f_1$ is once again a surjection. But now this new surjection takes $x_0$ into $\varphi_0$ and all remaining terms $t$ where $\psi$ acts into $\varphi_1$. Nevertheless, {\sc F11} grants us the existence of an injection $h_1$ whose $\mathfrak{F}$-domain is a restriction of $\psi$, whose $\mathfrak{F}$-image is $\varphi_2$, but such that $h_1$ {\em is not\/} a restriction of $f_1$. So, now we have $h_1(y_0) = \varphi_0$ and $h_1(y_1) = \varphi_1$, where $y_0$ and $y_1$ are terms where $\psi$ acts (and $\varphi_0$ and $\varphi_1$ are terms where $\varphi_2$ acts) but with no guarantee that either one of them is $x_0$. Actually, in this second step, we do have the guarantee that $h_1(x_0)\neq\varphi_0$ ($h_1(x_0)$ could be even $\underline{\mathfrak{0}}$ for all that matters). Function $g_1$, the inverse of $h_1$, is the injection demanded by the Partition Principle in Theorem \ref{teoremadeparticao}. Consequently, $g_0$ is no restriction of $g_1$. Thus, we have failed in fixing $x_0$ (associating it to $\varphi_0$) in order to fix a next term where $\psi$ acts and associate it to the next ordinal $\varphi_1$. At each step of our induction we have new terms associated to $\varphi_0$, $\varphi_1$, and so on, with the guarantee that at least one of them is not fixed to the same ordinal in the previous steps. Suppose now we reach the first limit ordinal $\omega$ in our already failed transfinite induction. Then we have a surjection $f_{\omega}:\psi\to\sigma_{\omega}$ such that $f_{\omega}(z_0) = \varphi_0$, $f_{\omega}(z_1) = \varphi_1$, and so on, and for all the remaining terms $t$ where $\psi$ acts we have $f_{\omega}(t) = \omega$, where $z_0$, $z_1$, and so on, are pairwise distinct terms where $\psi$ acts, but with no relationship whatsoever with the terms chosen by {\sc F11} in the previous steps for ordinals $o$ such that $\omega[o]$ (those $o$ who are strictly lesser than $\omega$). In the next immediate step of our induction a whole new plethora of terms $w_0$, $w_1$, $\cdots$, $w_{\omega}$, will be chosen by {\sc F11}, meaning we have no halting criterium to stop induction for the case when all terms of $\psi$ are exhausted. So, it is not possible to prove there is an injection between $\psi$ and an ordinal $o$ by transfinite induction, in the sense of fixing an injection $f:\psi\to o$ for some ordinal $o$ (through a process of transfinite recursion that is supposed to define $f$), specially when we already know that a hyperfunction can be `bigger' than $\omega$ (Theorem \ref{teoremaqueantecedegrande}). Finally, since we have no criterium for establishing a surjection with $\mathfrak{F}$-domain $o$ and $\mathfrak{F}$-image $\psi$ (in order to try another way of using {\sc F11}; but even if we had such a criterium, we would still face the same problems pointed above), that concludes our proof.
\end{description}

A similar proof of last theorem can be carried out by assuming any uncountable set like, e.g., the set of real numbers instead of a hyperfunction $\psi$.

If {\sc F11} demanded the chosen function $c$ {\em to be a restriction\/} of the surjection $f$, our {\em new axiom\/} of $\mathfrak{F}$-Choice {\sc F11'} could be easily used to prove any ZF-set (including hyperfunctions) can be well-ordered, and such a proof can be done by transfinite induction, in a similar way to Zermelo's Theorem. The idea is as follows. Let $x$ be a ZF-set and $f_0:x\to\varphi_1$ be given by $\forall t(x[t]\Rightarrow f_0(t) = \varphi_0)$. Theorem \ref{teoremadeparticao} (adapted to the new {\sc F11'}) grants a $g_0$ such that $g(x_0) = \varphi_0$, where $x_0$ is a term such that $x[x_0]$. For the next immediate step we define a function $f_1:x\to \varphi_2$ as $\forall t((f_0[t]\Rightarrow f_1(t) = f_0(t))\wedge ((x[t]\wedge \neg f_0[t])\Rightarrow f_1(t) = \varphi_1))$. Theorem \ref{teoremadeparticao} (adapted to the new {\sc F11'}) grants a $g_1$ such that $g_1(x_0) = \varphi_0$ and $g_1(x_1) = \varphi_1$, since our new proposed chosen function $c$ is a restriction of $f_1$. For the remaining steps, let $\kappa$ be an ordinal (Definition \ref{definicaodeordinais}). Then we can define a function $f_{\kappa}:x\to \sigma_{\kappa}$ such that $\forall \lambda \forall t(\kappa[\lambda]\Rightarrow (f_{\lambda}[t]\Rightarrow f_{\kappa}(t) = f_{\lambda}(t))\wedge ((x[t]\wedge\neg f_{\lambda}[t])\Rightarrow f_{\kappa}(t) = \kappa))$. Theorem \ref{teoremadeparticao} (adapted to the new {\sc F11'}) grants the existence of a $g_{\kappa}$ such that $g_{\kappa}(x_n) = n$ for any ordinal $n$ such that $\sigma_{\kappa}[n]$. We repeat the process until all terms where $x$ acts are exhausted. That is a recursive definition for a well order on $x$.

Thus, the issue here is the condition in {\sc F11} that the choice function $c$ is no restriction of $f$. {\sc F11} is strong enough to grant PP, but too weak to entail AC.

\section{Final remarks}\label{finalremarks}

As a reference to Heraclitus's flux doctrine, we are inclined to refer to all terms of Flow as {\em fluents\/}, rather than functions. In this sense Flow is a {\em theory of fluents\/}. That is also an auspicious homage to the {\em Method of Fluxions\/} by Isaac Newton \cite{Newton-64}. The famous `natural philosopher' referred to functions as fluents, and their derivatives as {\em fluxions\/}. Whether Newton was inspired or not by Heraclitus, that is historically uncertain (\cite{Royal-46}, page 38). Notwithstanding, we find such a coincidence quite inspiring and utterly opportune.

From the mathematical point of view, our framework was strongly motivated by von Neumann's set theory \cite{vonNeumann-25}, as we briefly discussed in the Introduction. But what are the main differences between Flow and von Neumann's ideas?

We refer to von Neumann's set theory as $\mbox{\boldmath{$\mathfrak{N}$}}$. The first important difference lurks in the way how von Neumann seemed to understand functions: ``a function can be regarded as a set of pairs, and a set as a function that can take two values... the two notions are completely equivalent''. That philosophical viewpoint seems to be committed to a set-theoretic framework. Here we follow a different path: {\sc i}) In $\mbox{\boldmath{$\mathfrak{N}$}}$ there are two privileged objects termed $A$ and $B$ which resemble our terms $\underline{\mathfrak{0}}$ and $\underline{\mathfrak{1}}$. Nevertheless, in $\mbox{\boldmath{$\mathfrak{N}$}}$ there is no further information about $A$ and $B$. Besides, those constants play a different role of $\underline{\mathfrak{0}}$ and $\underline{\mathfrak{1}}$ in $\mbox{\boldmath{$\mathfrak{F}$}}$, as we can see in the next item. {\sc ii}) In $\mbox{\boldmath{$\mathfrak{N}$}}$ there are two sorts of objects, namely, arguments and functions. Eventually some of those objects are both of them. And when an object $f$ is an argument and a function, which takes only values $A$ and $B$, then $f$ is a set. Our terms $\underline{\mathfrak{0}}$ and $\underline{\mathfrak{1}}$ have no similar role. Besides, we do not need to distinguish functions from any other kind of term. All objects of $\mbox{\boldmath{$\mathfrak{F}$}}$ are functions. {\sc iii}) Von Neumann believed a distinction between arguments and functions was necessary to avoid the well known antinomies from naive set theory. Nevertheless, we proved a simple axiom of self-reference ({\sc F2}) is enough to avoid such a problem. {\sc iv}) The distinction between sets and other collections which are `too big' to be sets depends on considerations if a specific term in $\mbox{\boldmath{$\mathfrak{N}$}}$ is both an argument and a function. Within $\mbox{\boldmath{$\mathfrak{F}$}}$ that distinction depends on considerations regarding $\mathfrak{F}$-successor. {\sc v}) Axiom I4 of $\mbox{\boldmath{$\mathfrak{N}$}}$ (\cite{vonNeumann-25}, page 399) says any function can be identified by its images. Our weak extensionality {\sc F1} allows us to derive a similar result as a non-trivial theorem (Theorem \ref{igualdadefuncoes}). {\sc vi}) One primitive concept in $\mbox{\boldmath{$\mathfrak{N}$}}$ is a binary functional letter (using modern terminology) which allows to define ordered pairs. In $\mbox{\boldmath{$\mathfrak{F}$}}$ that assumption is unnecessary. {\sc vii}) In $\mbox{\boldmath{$\mathfrak{N}$}}$ there are many constant functions (\cite{vonNeumann-25}, page 399, axiom II2). In $\mbox{\boldmath{$\mathfrak{F}$}}$ there are only two constant functions, namely, $\underline{\mathfrak{0}}$ and $\varphi_0$. {\sc viii}) In $\mbox{\boldmath{$\mathfrak{N}$}}$ the well-ordering theorem and the axiom of choice are consequences of its postulates. In $\mbox{\boldmath{$\mathfrak{F}$}}$ such a phenomenon does not take place.

Although we did not use $\mathfrak{F}$-composition to obtain our main result, we introduced that concept to deliver a partial response to the unavoidable question about how to `compose' functions in Flow. But other proposals for `composition' may be introduced. A simple example refers to what we call $\mathbb{Z}$-composition: $f$ and $g$ can be $\mathbb{Z}$-composed iff $Dom_f^{\mathfrak{F}}$, $Dom_g^{\mathfrak{F}}$, $Im_f^{\mathfrak{F}}$, and $Im_g^{\mathfrak{F}}$ act only on ZF-sets, and $Dom_g^{\mathfrak{F}} = Im_f^{\mathfrak{F}}$, and $f$ and $g$ do not act on each other. The $\mathbb{Z}$-composition $h = g\circ_{\mathbb{Z}} f$ is defined as $h(t) = g(f(t))$ for any $t$ where $f$ acts. It is easy to show $\mathbb{Z}$-composition is associative. In a sense, $\mathfrak{F}$-composition generalizes $\mathbb{Z}$-composition. We intend to study other possibilities for `composition' in forthcoming papers.

Other open problems refer to how precisely Flow is related to Lambda Calculus \cite{Barendregt-13}, Category Theory \cite{Hatcher-68} \cite{Lawvere-03}, String Diagrams \cite{Dixon-13}, and Autocategories \cite{Guitart-14}.

\section{Acknowledgements}

We thank Aline Zanardini, Bruno Victor, and Cl\'eber Barreto for insightful discussions which inspired this work. We acknowledge with thanks as well Newton da Costa, Jean-Pierre Marquis, Ed\'elcio Gon\c calves de Souza, D\'ecio Krause, Jonas Arenhart, Kherian Gracher, Bryan Leal Andrade, and B\'arbara Guerreira for valuable remarks and criticisms. We were benefitted also by comments and criticisms made by several participants of seminars delivered at Federal University of Paran\'a and Federal University of Santa Catarina. Finally we acknowledge Pedro D. Dam\'azio's provocation, more than 20 years ago, which worked as the first motivation towards the present work.

\end{document}